\documentclass[reqno]{amsart}
\usepackage{amssymb,latexsym}
\usepackage{amsfonts,mathrsfs}
\usepackage[margin=1.4in]{geometry}

\newtheorem{theorem}{Theorem}[section]
\newtheorem{lemma}[theorem]{Lemma}
\newtheorem{proposition}[theorem]{Proposition}
\newtheorem{corollary}[theorem]{Corollary}

\newcommand{\ba}{\begin{array}}
\newcommand{\ea}{\end{array}}

\def \qed{\cqfd}

\newcommand{\Vol}{\mathrm{Vol}}
\newcommand{\rank}{\mathrm{rank}}
\newcommand{\Id}{\mathrm{Id}}
\newcommand{\End}{\mathrm{End}}

\def\qed{\vbox{\hrule
\hbox{\vrule\hbox to 5pt{\vbox to 8pt{\vfil}\hfil}\vrule}\hrule}}

\newcommand{\beg}{\begin{eqnarray*}}
\newcommand{\begn}{\begin{eqnarray}}
\newcommand{\en}{\end{eqnarray*}}
\newcommand{\enn}{\end{eqnarray}}

\newcommand{\tr}{\mbox{\rm tr\,}}

\begin{document}
\title{Semi-stable Higgs sheaves and Bogomolov type inequality }
\author{Jiayu Li}
\address{School of Mathematical Sciences\\
University of Science and Technology of China\\
Hefei, 230026\\ and AMSS, CAS, Beijing, 100080, P.R. China\\} \email{jiayuli@ustc.edu.cn}
\author{Chuanjing Zhang}
\address{School of Mathematical Sciences\\
University of Science and Technology of China\\
Hefei, 230026,P.R. China\\ } \email{chjzhang@mail.ustc.edu.cn}
\author{Xi Zhang}
\address{School of Mathematical Sciences\\
University of Science and Technology of China\\
Hefei, 230026,P.R. China\\ } \email{mathzx@ustc.edu.cn}
\thanks{The authors were supported in part by NSF in
China,  No. 11571332, 11131007, 11526212.}

\begin{abstract} In this paper, we study semistable Higgs sheaves over compact K\"ahler manifolds, we prove that there is an approximate admissible Hermitian-Einstein structure
on a semi-stable reflexive Higgs sheaf and consequently, the Bogomolove type inequality holds on a semi-stable reflexive Higgs sheaf.
\end{abstract}

\maketitle

\section{Introduction}
\setcounter{equation}{0}

\hspace{0.4cm}

Let $(M, \omega )$ be a compact K\"ahler manifold, and $E$ be a holomorphic vector bundle on $M$. Donaldson-Uhlenbeck-Yau theorem states that the $\omega$-stability of $E$ implies the existence of  $\omega$-Hermitian-Einstein metric on $E$. Hitchin \cite{Hi} and Simpson \cite{Si} proved that the theorem holds also for Higgs bundles. We \cite{LZ2} proved that there is an approximate  Hermitian-Einstein structure on a semi-stable Higgs bundle, which confirms a conjecture due to Kobayashi \cite{Ko2} (also see \cite{Ja}). There are many interesting and important works related
(\cite{LY, Hi, Si, Br1, BS, GP, BG, AG1, Bi, BT, LN1, LN2, M, Mo1, Mo}, etc.). Among all of them, we recall that, Bando and Siu \cite{BS} introduced the notion of admissible Hermitian metrics on torsion-free sheaves, and proved the Donaldson-Uhlenbeck-Yau theorem on stable reflexive sheaves.

Let $\mathcal{E}$ be a  torsion-free coherent sheaf, and $\Sigma $ be the set of singularities where $\mathcal{E}$ is not locally free. A Hermitian metric $H$ on the holomorphic bundle $\mathcal{E}|_{M\setminus \Sigma }$ is called {\it admissible } if

(1) $|F_{H}|_{H, \omega}$ is square integrable;

(2) $|\Lambda_{\omega } F_{H}|_{H}$ is uniformly bounded.\\
Here $F_{H}$ is the curvature tensor of Chern connection $D_{H}$ with respect to the Hermitian metric $H$, and $\Lambda_{\omega }$ denotes the contraction with  the K\"ahler metric $\omega $.

   Higgs bundle and Higgs sheaf are studied by Hitchin (\cite{Hi}) and Simpson (\cite{Si}, \cite{Si2}), which play an important role in many different areas including gauge theory, K\"ahler and hyperk\"ahler geometry, group representations, and nonabelian Hodge theory.
A Higgs sheaf on $(M, \omega )$ is a pair $(\mathcal{E}, \phi )$ where $\mathcal{E}$ is a  coherent sheaf  on $M$ and the Higgs field $\phi \in \Omega^{1,0 }(\End (\mathcal{E}))$ is a holomorphic section such that $\phi \wedge \phi =0 $. If the sheaf $\mathcal{E}$ is torsion-free (resp. reflexive, locally free), then we say the Higgs sheaf $(\mathcal{E}, \phi )$ is torsion-free (resp. reflexive, locally free). A torsion-free Higgs sheaf $(\mathcal{E}, \phi )$ is said to be
  $\omega $-stable (respectively, $\omega$-semi-stable), if for every $\phi$-invariant coherent proper sub-sheaf
$\mathcal{F}\hookrightarrow \mathcal{E}$, it holds:
\begin{eqnarray}
\mu_{\omega} (\mathcal{F})=\frac{\deg_{\omega} (\mathcal{F})}{\rank (\mathcal{F})}< (\leq ) \mu_{\omega} (\mathcal{E})=\frac{\deg_{\omega} (\mathcal{E})}{\rank (\mathcal{E})},
\end{eqnarray}
where $\mu_{\omega} (\mathcal{F})$ is called the $\omega$-slope of $\mathcal{F}$.

Given a Hermitian metric $H$ on the locally free part of the  Higgs sheaf $(\mathcal{E}, \phi )$, we consider the  Hitchin-Simpson connection
 \begin{equation}\overline{\partial}_{\phi}:=\overline{\partial}_{\mathcal{E}}+\phi , \quad D_{H,  \phi }^{1, 0}:=D_{H }^{1, 0}  +\phi ^{\ast H}, \quad D_{H,  \phi }= \overline{\partial}_{\phi}+ D_{H,  \phi }^{1, 0},\end{equation}
where $D_{H}$ is the Chern connection with respect to the metric $H$ and $\phi ^{\ast H}$ is the adjoint of $\phi $ with respect to  $H$. The curvature of the  Hitchin-Simpson connection is
\begin{equation}
F_{H, \phi}=F_{H} +[\phi , \phi ^{\ast H}] +D_{H}^{1, 0}\phi + \overline{\partial }_{\mathcal{E}} \phi^{\ast H},
\end{equation}
where $F_{H}$ is the curvature of the Chern connection $D_{H}$. A Hermitian metric $H$ on the  Higgs sheaf $(\mathcal{E}, \phi )$ is said to be
admissible Hermitian-Einstein  if it is admissible and
  satisfies the following Einstein condition on $M\setminus \Sigma $, i.e
\begin{equation}\label{HE}
\sqrt{-1}\Lambda_{\omega} (F_{H} +[\phi , \phi ^{\ast H}])
=\lambda \Id_{\mathcal{E}},
\end{equation}
where   $\lambda $ is a constant
given by $\lambda  =\frac{2\pi}{\Vol(M, \omega)} \mu_{\omega} (\mathcal{E})$.
Hitchin (\cite{Hi}) and Simpson (\cite{Si}) proved that a Higgs bundle admits
a Hermitian-Einstein metric if and only if it's Higgs poly-stable. Biswas and Schumacher \cite{BiS} studied the Donaldson-Uhlenbeck-Yau theorem for reflexive Higgs sheaves.

In this paper, we study the semi-stable Higgs sheaves. We say a torsion-free  Higgs sheaf $(\mathcal{E}, \phi )$ admits an approximate admissible Hermitian-Einstein structure if for every positive $\delta
$, there is an admissible Hermitian metric $H_{\delta }$ such that
\begin{equation}
\sup _{x\in M\setminus \Sigma } |\sqrt{-1}\Lambda_{\omega }(F_{H_{\delta}}+[\phi , \phi^{\ast H_{\delta}}])-\lambda \Id_{\mathcal{E}}|_{H_{\delta}}(x)<\delta .
\end{equation}
The approximate Hermitian-Einstein structure was introduced by Kobayashi (\cite{Ko2}) on a holomorphic vector bundle, it is the differential geometric counterpart of the semi-stability. Kobayashi \cite{Ko2} proved there is an approximate Hermitian-Einstein structure on a semi-stable holomorphic vector bundle over an algebraic manifold, which he conjectured should be true over any K\"ahler manifold. The conjecture was confirmed in \cite{Ja,LZ2}.  In this paper, we proved our theorem holds for a semi-stable reflexive Higgs sheaf over a compact K\"ahler manifold.

\medskip

\begin{theorem}\label{thm 1.1}
A reflexive Higgs sheaf $(\mathcal{E}, \phi )$  on an $n$-dimensional compact K\"ahler manifold $(M, \omega)$ is semi-stable, if and only if it admits an approximate admissible Hermitian-Einstein structure. Specially, for a semi-stable reflexive Higgs sheaf $(\mathcal{E}, \phi )$ of rank $r$, we have the following Bogomolov type inequality
\begin{equation}\label{Bog1}
\int_{M} (2c_{2}(\mathcal{E})-\frac{r-1}{r}c_{1}(\mathcal{E})\wedge c_{1}(\mathcal{E}))\wedge\frac{\omega^{n-2}}{(n-2)!}\geq 0 .
\end{equation}
\end{theorem}
\medskip

The Bogomolov inequality was first obtained by Bogomolov (\cite{Bo}) for semi-stable holomorphic vector bundles over complex algebraic surfaces, it had been extended to certain
classes of generalized vector bundles, including parabolic bundles and orbibundles. By constructing a Hermitian-Einstein metric, Simpson proved the Bogomolov inequality for stable Higgs bundles on compact K\"ahler
manifolds. Recently, Langer (\cite{La}) proved the Bogomolov type inequality for semi-stable Higgs sheaves over algebraic varieties by using an algebraic-geometric method. His  method can not be applied to the K\"ahler manifold case. We use analytic method to study the Bogomolov inequality for semi-stable reflexive Higgs sheaves over compact K\"ahle manifolds, new idea is needed.

\medskip

  We now give an overview of our proof. As in \cite{BS}, we make a regularization on the reflexive sheaf $\mathcal{E}$, i.e. take blowing up with smooth centers finite times $\pi_{i}: M_{i}\rightarrow M_{i-1}$, where $i=1, \cdots , k$ and $M_{0}=M$, such that the pull-back of $\mathcal{E}^{\ast}$ to $M_{k} $ modulo torsion  is locally free and
  \begin{equation}
  \pi = \pi_{1}\circ \cdots \circ \pi_{k}: M_{k}\rightarrow M
  \end{equation}
  is biholomorphic outside $\Sigma $. In the following, we denote $M_{k}$ by $\tilde{M}$,  the exceptional divisor $\pi^{-1} \Sigma $ by $\tilde{\Sigma }$, and the holomorphic vector bundle $(\pi^{\ast}\mathcal{E}^{\ast}/torsion)^{\ast}$ by $E$. Since $\mathcal{E}$ is locally free outside $\Sigma $, and the holomorphic bundle $E$ is isomorphic to $\mathcal{E}$ on $\tilde{M}\setminus \tilde{\Sigma}$,  the pull-back field $\pi^{\ast}\phi $ is a holomorphic section of $\Omega^{1,0 }(\End (E))$ on $\tilde{M}\setminus \tilde{\Sigma}$. By Hartogs' extension theorem, the holomorphic section $\pi^{\ast}\phi $ can be extended to the whole $\tilde{M}$ as a Higgs field of $E$. In the following, we also denote the extended Higgs field $\pi^{\ast}\phi $ by $\phi $ for simplicity. So we get a Higgs bundle $(E, \phi )$ on $\tilde{M}$ which is isomorphic to the Higgs sheaf $(\mathcal{E}, \phi )$ outside the exceptional divisor $\tilde{\Sigma }$.

   It is well known that $\tilde{M}$ is also K\"ahler (\cite{GH}). Fix a K\"ahler metric $\eta $ on $\tilde{M}$ and set
 \begin{equation}\omega_{\epsilon }=\pi^{\ast}\omega +\epsilon \eta \end{equation}
  for any small $0<\epsilon \leq 1$. Let $K_{\epsilon}(t, x, y)$ be the heat kernel with respect to the K\"ahler metric $\omega_{\epsilon}$. Bando and Siu (Lemma 3 in \cite{BS}) obtained a uniform  Sobolev inequality for $(\tilde{M}, \omega_{\epsilon})$, using Cheng and Li's estimate (\cite{CL}), they got a uniform upper bound of the heat kernels $K_{\epsilon}(t, x, y)$.
 Given a smooth Hermitian metric $\hat{H}$ on the bundle $E$, it is easy to see that there exists a constant $\hat{C}_{0}$ such that
 \begin{equation}\label{initial1}
 \int_{\tilde{M}}( |\Lambda_{\omega_{\epsilon }}F_{\hat{H}}|_{\hat{H}}+|\phi |_{\hat{H}, \omega_{\epsilon}}^{2})\frac{\omega_{\epsilon}^{n}}{n!}\leq \hat{C}_{0},
 \end{equation}
 for all $0< \epsilon \leq 1$. This also gives a uniform bound on $ \int_{\tilde{M}} |\Lambda_{\omega_{\epsilon }}(F_{\hat{H}}+[\phi , \phi ^{\ast \hat{H}}])|_{\hat{H}}\frac{\omega_{\epsilon}^{n}}{n!}$.

We study the following evolution equation on Higgs bundle $(E, \phi )$ with the fixed initial metric $\hat{H}$ and with respect to the K\"ahler metric $\omega_{\epsilon}$,
\begin{equation}\label{DDD1}
\left \{\begin{split} &H_{\epsilon}(t)^{-1}\frac{\partial H_{\epsilon}(t)}{\partial
t}=-2(\sqrt{-1}\Lambda_{\omega_{\epsilon}}(F_{H_{\epsilon}(t)}+[\phi , \phi ^{\ast H_{\epsilon}(t)}])-\lambda_{\epsilon }\Id_{E}),\\
&H_{\epsilon}(0)=\hat{H},\\
\end{split}
\right.
\end{equation}
where  $\lambda_{\epsilon}  =\frac{2\pi}{\Vol(\tilde{M}, \omega_{\epsilon})} \mu_{\omega_{\epsilon}} (E)$. Simpson (\cite{Si}) proved the existence of long time solution of the above heat flow.  By the standard parabolic estimates and  the uniform upper bound of the heat kernels $K_{\epsilon}(t, x, y)$, we know that $|\Lambda_{\omega_{\epsilon }}(F_{H_{\epsilon}(t)}+[\phi , \phi ^{\ast H_{\epsilon}(t)}])|_{H_{\epsilon}(t)}$ has a uniform $L^{1}$ bound for $t\geq 0$ and a uniform $L^{\infty}$ bound for $t\geq t_{0}>0$. As in \cite{BS}, taking the limit as $\epsilon \rightarrow 0$, we have a long time solution $H(t)$ of the following evolution equation on $M\setminus \Sigma \times [0, +\infty)$, i.e. $H(t)$ satisfies:
\begin{equation}\label{SSS1}
\left \{\begin{split} &H(t)^{-1}\frac{\partial H(t)}{\partial
t}=-2(\sqrt{-1}\Lambda_{\omega}(F_{H(t)}+[\phi , \phi ^{\ast H(t)}])-\lambda \Id_{\mathcal{E}}),\\
&H(0)=\hat{H}.\\
\end{split}
\right.
\end{equation}
Here $H(t)$ can be seen as a Hermitian metric  defined on the locally free part of $\mathcal{E}$, i.e. on $M\setminus \Sigma$.

In order to get the admissibility of Hermitian metric $H(t)$ for positive time $t>0$, we should show that $|\phi|_{H(t), \omega }\in L^{\infty}$ for $t>0$. In fact, we can prove that  $|\phi|_{H(t), \omega }$ has a uniform $L^{\infty}$ bound for $t\geq t_{0}>0$. In \cite{LZ1}, by using the maximum principle, we proved this uniform $L^{\infty}$ bound of $|\phi|_{H(t), \omega }$ along the evolution equation for the Higgs bundle case. In the Higgs sheaf case, since the equation (\ref{SSS1}) has singularity on $\Sigma $, we can not use the maximum principle directly. So we need new argument to get a uniform $L^{\infty}$ bound of $|\phi|_{H(t), \omega }$, see section 3 for details.

 The key part in the proof of Theorem \ref{thm 1.1} is to prove the existence of admissible approximate Hermitian-Einstein structure on a semi-stable reflexive Higgs sheaf. The Bogomolov type inequality (\ref{Bog1}) is an application. In fact, we prove that if the reflexive Higgs sheaf $(\mathcal{E}, \phi )$ is semi-stable, along the evolution equation (\ref{SSS1}), we must have
\begin{equation}\label{SS01}
\sup _{x\in M\setminus \Sigma } |\sqrt{-1}\Lambda_{\omega }(F_{H(t)}+[\phi , \phi^{\ast H(t)}])-\lambda \Id_{\mathcal{E}}|_{H(t)}(x)\rightarrow 0,
\end{equation}
as $t\rightarrow +\infty $.  We prove (\ref{SS01}) by contradiction, if not,  we can construct a saturated Higgs subsheaf such that its $\omega $-slope is greater than $\mu_{\epsilon}(\mathcal{E})$. Since the singularity set $\Sigma $ is a complex analytic subset with co-dimension at least $3$, it is easy to show that $(M\setminus \Sigma , \omega )$ satisfies all three assumptions that Simpson (\cite{Si}) imposes on the non-compact base K\"ahler manifold. Let's recall Simpson's argument for a Higgs bundle in the case where the base K\"ahler manifold is non-compact. Simpson assumes that there exists a good initial  Hermitian metric $K$ satisfying $\sup_{M\setminus \Sigma}|\Lambda_{\omega}F_{K, \phi }|_{K}<\infty $, then he defines the analytic stability for $(\mathcal{E}, \phi , K)$ by using the Chern-Weil formula with respect to the metric $K$ (Lemma 3.2 in \cite{Si}). Under the $K$-analytic stability condition, he constructs a Hermitian-Einstein metric for the Higgs bundle by limiting the evolution equation (\ref{SSS1}).

Here, we have to pay more attention to the analytic stability (or semi-stability) of $(\mathcal{E}, \phi )$.
 Let $\mathcal{F}$ be a saturated sub-sheaf of $\mathcal{E}$,   we know that $\mathcal{F}$ can be seen as a sub-bundle of $\mathcal{E}$ outside a singularity set $V=\Sigma_{\mathcal{F}}\cup \Sigma $ of codimension at least $2$, then $\hat{H}$ induces a Hermitian metric $\hat{H}_{\mathcal{F}}$ on $\mathcal{F}$. Bruasse (Proposition 4.1 in \cite{Br}) had proved the following Chern-Weil formula
 \begin{eqnarray}\label{CW1}
 \deg_{\omega }(\mathcal{F})=\int_{M\setminus V} c_{1}(\mathcal{F}, \hat{H}_{\mathcal{F}})\wedge \frac{\omega^{n-1}}{(n-1)!},
 \end{eqnarray}
where $c_{1}(\mathcal{F}, \hat{H}_{\mathcal{F}})$ is the first Chern form with respect to the induced metric $\hat{H}_{\mathcal{F}}$.  By (\ref{CW1}), we see that the stability (semi-stability ) of the reflexive Higgs sheaf $(\mathcal{E}, \phi )$ is  equivalent to the analytic stability (semi-stability) with respect to the metric $\hat{H}$ in Simpson's sense. But, we are not clear whether the above Chern-Weil formula is still valid if the metric $\hat{H}$ is replaced by an admissible metric $H(t)$ ($t>0$). So, the stability (or semi-stability) of the reflexive Higgs sheaf $(\mathcal{E}, \phi )$ may not imply the analytic stability (or semi-stability ) with respect to the metric $H(t)$ ($t>0$).   The admissible metric $H(t)$ ($t>0$) can not be chosen as a good initial metric in Simpson's sense.  On the other hand,  the initial metric $\hat{H}$ may not satisfy the curvature finiteness condition (i.e. $|\Lambda_{\omega}F_{\hat{H}, \phi }|_{\hat{H}}$ may not be $L^{\infty}$ bounded), so we should modify Simpson's argument in our case, see the proof of Proposition \ref{prop 4.1} in section 4 for details.

\medskip

  If the reflexive Higgs sheaf $(\mathcal{E} , \phi )$ is $\omega $-stable, it is well known that the pulling back Higgs bundle $(E, \phi )$ is $\omega_{\epsilon }$-stable for sufficiently small $\epsilon$. By Simpson's result (\cite{Si}), there exists an $\omega_{\epsilon }$-Hermitian-Einstein metric $H_{\epsilon}$ for every small $\epsilon$. In \cite{BS}, Bando and Siu point out that it is possible to get an $\omega$-Hermitian-Einstein metric $H$ on the reflexive Higgs sheaf $(\mathcal{E} , \phi )$ as a limit of  $\omega_{\epsilon }$-Hermitian-Einstein metric $H_{\epsilon}$ of Higgs bundle $(E, \phi )$ on $\tilde{M}$ as $\epsilon \rightarrow 0$. In the end of this paper, we solve this problem.

  \medskip

\begin{theorem}\label{thm 1.2}
Let $H_{\epsilon}$ be an $\omega_{\epsilon }$-Hermitian-Einstein metric on the Higgs bundle $(E, \phi )$,  by choosing a subsequence and rescaling it, $H_{\epsilon}$ must converge to an $\omega$-Hermitian-Einstein metric $H$ in local $C^{\infty}$-topology outside the exceptional divisor $\tilde{\Sigma }$ as $\epsilon \rightarrow 0$.
\end{theorem}
  \medskip

This paper is organized as follows. In Section 2, we recall some basic estimates for the heat flow (\ref{DDD1}) and give proofs for local uniform $C^{0}$, $C^{1}$ and higher order estimates for reader's convenience.  In section 3, we give a uniform $L^{\infty}$ bound for the norm of the Higgs field along the heat flow (\ref{SSS1}).   In section 4, we prove the existence of admissible approximate Hermitian-Einstein structure on the semi-stable reflexive Higgs sheaf  and complete the proof of Theorem \ref{thm 1.1}. In section 5,  we prove Theorem \ref{thm 1.2}.

\hspace{0.3cm}

\section{Analytic preliminaries and basic estimates }
\setcounter{equation}{0}

 Let $(M, \omega )$ be a compact K\"ahler manifold of complex dimension $n$, and $(\mathcal{E}, \phi )$  be a reflexive Higgs
sheaf on $M$ with the singularity set $\Sigma $. There exists a bow-up $\pi : \tilde{M}\rightarrow M$ such that the pulling back Higgs bundle $(E, \phi )$ on $\tilde{M}$  is isomorphic to $(\mathcal{E}, \phi )$ outside the exceptional divisor $\tilde{\Sigma }=\pi^{-1}\Sigma $.  It is well known that $\tilde{M}$ is also K\"ahler (\cite{GH}). Fix a K\"ahler metric $\eta $ on $\tilde{M}$ and set $\omega_{\epsilon }=\pi^{\ast}\omega +\epsilon \eta $
  for $0<\epsilon \leq 1$. Let $K_{\epsilon}(x, y, t)$ be the heat kernel with respect to the K\"ahler metric $\omega_{\epsilon}$. Bando and Siu (Lemma 3 in \cite{BS}) obtained a uniform  Sobolev inequality for $(\tilde{M}, \omega_{\epsilon})$. Combining Cheng and Li's estimate (\cite{CL}) with  Grigor'yan's result (Theorem 1.1 in \cite{Gr}), we have the following uniform upper bound of the heat kernels, furthermore, we also have  a uniform lower bound of the Green functions.

 \medskip

\begin{proposition}\label{Prop 2.1}{\bf(Proposition 2 in \cite{BS})}
Let $K_{\epsilon}$ be the heat kernel with respect to the metric $\omega_{\epsilon}$, then for any $\tau >0$, there exists a constant $C_{K}(\tau)$ which is independent of $\epsilon $, such that
\begin{equation}\label{kernel01}0\leq K_{\epsilon}(x , y, t)\leq C_{K}(\tau) (t^{-n}\exp{(-\frac{(d_{\omega_{\epsilon}}(x, y))^{2}}{(4+\tau )t})}+1)\end{equation} for every $x, y\in \tilde{M}$ and $0<t < +\infty$, where $d_{\omega_{\epsilon}}(x, y)$ is the distance between $x$ and $y$ with respect to the metric $\omega_{\epsilon}$. There also exists a constant $C_{G}$ such that
\begin{equation}\label{green} G_{\epsilon}(x , y)\geq -C_{G}\end{equation}
for every $x, y\in \tilde{M}$ and $0<\epsilon \leq 1$, where $G_{\epsilon}$ is the Green function with respect to the metric $\omega_{\epsilon}$.
\end{proposition}
 \medskip

Let $H_{\epsilon } (t)$  be the long time solutions of the heat flow (\ref{DDD1}) on the Higgs bundle $(E, \phi )$ with the fixed smooth initial metric $\hat{H}$ and with respect to the K\"ahler metric $\omega_{\epsilon}$. By (\ref{initial1}), there is a constant $\hat{C}_{1}$ independent of $\epsilon $ such that
\begin{equation}
\int_{\tilde{M}} |\sqrt{-1}\Lambda_{\omega_{\epsilon }}(F_{\hat{H}}+[\phi , \phi ^{\ast \hat{H}}])-\lambda_{\epsilon }\Id_{E}|_{\hat{H}}\frac{\omega_{\epsilon}^{n}}{n!}\leq \hat{C}_{1}.
\end{equation}
 For simplicity, we set:
\begin{equation}
\Phi (H_{\epsilon}(t), \omega_{\epsilon})=\sqrt{-1}\Lambda_{\omega_{\epsilon}}(F_{H_{\epsilon}(t)}+[\phi , \phi ^{\ast H_{\epsilon}(t)}])-\lambda_{\epsilon }\Id_{E} .
\end{equation}
The following estimates are essentially proved by Simpson (Lemma 6.1 in \cite{Si}, see also Lemma 4 in \cite{LZ2}).
Along  the heat flow (\ref{DDD1}), we have:
 \begin{equation}\label{F1}
(\Delta_{\epsilon}-\frac{\partial }{\partial t})\tr (\Phi (H_{\epsilon}(t), \omega_{\epsilon}))=0,
\end{equation}
\begin{equation}\label{F2}
(\Delta_{\epsilon}-\frac{\partial }{\partial t} )|\Phi (H_{\epsilon}(t), \omega_{\epsilon})|_{H_{\epsilon}(t)}^{2}=2|D_{H_{\epsilon}, \phi} (\Phi (H_{\epsilon}(t), \omega_{\epsilon}))|^{2}_{H_{\epsilon}(t), \omega_\epsilon},
\end{equation}
and
\begin{equation}\label{H00}
(\Delta_{\epsilon } -\frac{\partial }{\partial t}) |\Phi (H_{\epsilon}(t), \omega_{\epsilon})|_{H_{\epsilon}(t)}\geq  0.
\end{equation}
Then, for $t>0$,
\begin{equation}\label{H001}
\int_{\tilde{M}}|\Phi (H_{\epsilon}(t), \omega_{\epsilon}) |_{H_{\epsilon}(t)}\frac{\omega_{\epsilon }^{n}}{n!}\leq \int_{\tilde{M}}|\Phi (\hat{H}, \omega_{\epsilon}) |_{\hat{H}}\frac{\omega_{\epsilon }^{n}}{n!}\leq \hat{C}_{1},
\end{equation}
 \begin{equation}\label{H000}
 \max_{x\in \tilde{M}}|\Phi (H_{\epsilon}(t), \omega_{\epsilon}) |_{H_{\epsilon}(t)}(x) \leq \int_{\tilde{M}}K_{\epsilon} (x, y, t)|\Phi (\hat{H}, \omega_{\epsilon}) |_{\hat{H}}\frac{\omega_{\epsilon }^{n}}{n!} ,
\end{equation}
and
\begin{equation}\label{H0001}
 \max_{x\in \tilde{M}}|\Phi (H_{\epsilon}(t+1), \omega_{\epsilon}) |_{H_{\epsilon}(t+1)}(x) \leq \int_{\tilde{M}}K_{\epsilon} (x, y, 1)|\Phi (H_{\epsilon}(t), \omega_{\epsilon}) |_{H_{\epsilon}(t)}\frac{\omega_{\epsilon }^{n}}{n!} .
\end{equation}

By the upper bound of the heat kernels (\ref{kernel01}), we have
 \begin{equation}\label{H005}
 \max_{x\in \tilde{M}}|\Phi (H_{\epsilon}(t), \omega_{\epsilon}) |_{H_{\epsilon}(t)}(x) \leq C_{K}(\tau )\hat{C}_{1} (t^{-n}+1),
\end{equation}
and
\begin{equation}\label{H00001}
 \max_{x\in \tilde{M}}|\Phi (H_{\epsilon}(t+1), \omega_{\epsilon}) |_{H_{\epsilon}(t+1)}(x) \leq 2C_{K}(\tau )\int_{\tilde{M}}|\Phi (H_{\epsilon}(t), \omega_{\epsilon}) |_{H_{\epsilon}(t)}\frac{\omega_{\epsilon }^{n}}{n!} .
\end{equation}

Set \begin{equation}\exp (S_{\epsilon}(t))=h_{\epsilon}(t)=\hat{H}^{-1}H_{\epsilon }(t),\end{equation}
 where $S_{\epsilon}(t)\in \End(E)$  is self-adjoint with respect to $\hat{H} $ and $H_{\epsilon}(t)$.
 By the heat flow (\ref{DDD1}), we have:
  \begin{equation}\frac{\partial }{\partial t} \log \det(h_\epsilon(t))=\tr (h_\epsilon^{-1}\frac{\partial h_\epsilon}{\partial t})=-2\tr (\Phi (H_{\epsilon}(t), \omega_{\epsilon})),\end{equation}
and 
\begin{eqnarray}\label{1}
\int_{\tilde{M}}\tr (S_{\epsilon}(t))\frac{\omega_{\epsilon }^{n}}{n!} =\int_{\tilde{M}}\log \det (h_{\epsilon}(t))\frac{\omega_{\epsilon }^{n}}{n!}=0
\end{eqnarray}
for all $t\geq 0$.

In the following, we denote:
\begin{equation}
B_{\omega_{1}}(\delta )=\{x\in \tilde{M}|d_{\omega_{1}}(x, \Sigma )<\delta \},
\end{equation}
where $d_{\omega_{1}}$ is the distance function with respect to the K\"ahler metric $\omega_{1}$. Since $\hat{H}$ is a smooth Hermitian metric on $E$, $\phi \in \Omega^{1,0}_{\tilde{M}}(\End(E))$ is a smooth field, and $\pi^\ast\omega$ is degenerate only along $\Sigma$, there exist constants $\hat{c}(\delta^{-1})$ and $\hat{b}_k(\delta^{-1})$  such that
\begin{equation}\label{initial2}
\begin{split}
&\{|\Lambda_{\omega_{\epsilon}}F_{\hat{H}}|_{\hat{H}}+|\phi|_{\hat{H}, \omega_{\epsilon}}^{2}\}(y)\leq \hat{c}(\delta^{-1}), \\
&\{|\nabla_{\hat{H}}^kF_{\hat{H}}|_{\hat{H}, \omega_{\epsilon}}^{2}+|\nabla_{\hat{H}}^{k+1}\phi|_{\hat{H}, \omega_{\epsilon}}^{2}\}\leq \hat{b}_k(\delta^{-1}),\\
\end{split}
\end{equation}
for all $y\in \tilde{M}\setminus B_{\omega_1}(\frac{\delta}{2})$, all $0\leq \epsilon \leq 1$ and all $k\geq 0$.

In order to get a uniform local $C^0$-estimate of $h_{\epsilon}(t)$,
We first prove that $|\Phi (H_{\epsilon}(t), \omega_{\epsilon})|_{H_{\epsilon}(t)}$ is uniform locally bounded, i.e. we obtain the following Lemma.

\medskip

\begin{lemma}\label{lem 2.2} There exists a constant $\tilde{C}_{1}(\delta ^{-1})$  such that
\begin{equation}\label{220}
|\Phi (H_{\epsilon}(t), \omega_{\epsilon})|_{H_{\epsilon}(t)}(x)\leq \tilde{C}_{1}(\delta ^{-1})
\end{equation}
for all $(x, t)\in (\tilde{M}\setminus B_{\omega_1}(\delta))\times [0, \infty )$, and all $0< \epsilon \leq 1$.
\end{lemma}

\medskip

{\bf Proof. }Using the inequality (\ref{H000}), we have
\begin{equation}\label{221}
|\Phi (H_{\epsilon}(t), \omega_{\epsilon})|_{H_{\epsilon}(t)}(x)\leq \Big(\int_{M\setminus B_\epsilon(\frac{\delta}{2})}+\int_{B_\epsilon(\frac{\delta}{2})}\Big)K_\epsilon(x, y, t)|\Phi (\hat{H}, \omega_{\epsilon})|_{\hat{H}}(y)\frac{\omega_{\epsilon}^n(y)}{n!}.
\end{equation}
Noting $\int_{\tilde{M}}K_\epsilon(x, y, t)\frac{\omega_{\epsilon}^n}{n!}=1$ and using (\ref{initial2}), we have
\begin{equation}\label{222}
\begin{split}
&\int_{\tilde{M}\setminus B_\epsilon(\frac{\delta}{2})}K_\epsilon(x, y, t)|\Phi (\hat{H}, \omega_{\epsilon})|_{\hat{H}}(y)\frac{\omega_{\epsilon}^n}{n!}\\
\leq&\ (\hat{c}(\delta^{-1})+ \lambda_\epsilon\sqrt{r})\int_{\tilde{M}}K_\epsilon(x, y, t)\frac{\omega_{\epsilon}^n(y)}{n!}\\
\leq&\ \hat{c}_1(\delta^{-1}).\\
\end{split}
\end{equation}
where $\hat{c}_1(\delta^{-1})$ is a constant independent of $\epsilon$. Since $\pi^\ast\omega$ is degenerate only along $\Sigma$, there exists a constant $\tilde{a}(\delta)$ such that
\begin{equation}\label{metric02}
\tilde{a}(\delta)\omega_1< \pi^\ast\omega < \omega_\epsilon < \omega_1
\end{equation}
on $\tilde{M}\setminus B_{\omega_1}(\frac{\delta}{4})$, for all $0< \epsilon \leq 1$. Let $x \in \tilde{M}\setminus B_{\omega_1}(\delta)$ and $y \in \partial (B_{\omega_1}(\frac{\delta}{2}))$, it is clear that
\begin{equation}
d_{\omega_\epsilon}(x, y)\geq d_{\pi^\ast\omega}(x, y)> \sqrt{\tilde{a}(\delta)}d_{\omega_1}(x, y)\geq \frac{\delta\sqrt{\tilde{a}(\delta)}}{2}.
\end{equation}
Let $a(\delta)=\frac{\delta\sqrt{\tilde{a}(\delta)}}{2}$.  If $x\in \tilde{M}\setminus B_{\omega_1}(\delta)$ and $y\in B_{\omega_1}(\frac{\delta}{2})$, we have
\begin{equation}
d_{\omega_\epsilon}(x, y)\geq a(\delta)
\end{equation}
for all $0\leq \epsilon \leq 1$. Then,
\begin{equation}\label{223}
\begin{split}
&\int_{B_{\omega_1}(\frac{\delta}{2})}K_\epsilon(x, y, t)|\Phi (\hat{H}, \omega_{\epsilon})|_{\hat{H}}(y)\frac{\omega_{\epsilon}^n(y)}{n!}\\
\leq &\ C_k(\tau )\int_{B_{\omega_1}(\frac{\delta}{2})}(t^{-n}\exp(-\frac{d_{\omega_\epsilon}(x, y)}{(4+\tau)t})+1)|\Phi (\hat{H}, \omega_{\epsilon})|_{\hat{H}}(y)\frac{\omega_{\epsilon}^n(y)}{n!}\\
\leq &\ C_k(\tau )\int_{B_{\omega_1}(\frac{\delta}{2})}(t^{-n}\exp(-\frac{a(\delta)}{(4+\tau)t})+1)
|\Phi (\hat{H}, \omega_{\epsilon})|_{\hat{H}}\frac{\omega_{\epsilon}^n}{n!}\\
\leq &\ C_k(\tau )\big(\frac{a(\delta)}{4+\tau}n\big)^{-n}\exp(-n)\int_{B_{\omega_1}(\frac{\delta}{2})}
|\Phi (\hat{H}, \omega_{\epsilon})|_{\hat{H}}\frac{\omega_{\epsilon}^n}{n!}\\
\leq &\ C_k(\tau )\hat{C}_{1}\big(\frac{a(\delta)}{4+\tau}n\big)^{-n}\exp(-n),\\
\end{split}
\end{equation}
for all $(x, t)\in (\tilde{M}\setminus B_{\omega_1}(\delta))\times [0, \infty )$. It is obvious that (\ref{221}), (\ref{222}) and (\ref{223}) imply (\ref{220}).

\hfill $\Box$ \\

By a direct calculation, we have
\begin{equation}\label{c01}
\begin{split}
&\frac{\partial}{\partial t}\log(\tr h_\epsilon(t)+\tr h_\epsilon^{-1}(t))\\
=&\ \frac{\tr(h_\epsilon(t)\cdot h_\epsilon^{-1}(t)\frac{\partial h_\epsilon(t)}{\partial t})-\tr(h_\epsilon^{-1}(t)\frac{\partial h_\epsilon(t)}{\partial t}\cdot h_\epsilon^{-1}(t))}{\tr h_\epsilon(t)+\tr h_\epsilon^{-1}(t)}\\
\leq &\ 2|\Phi (H_{\epsilon}(t), \omega_{\epsilon})|_{H_\epsilon(t)},
\end{split}
\end{equation}
and
\begin{equation}\label{c02}
\log (\frac{1}{2r}(\tr h_\epsilon(t) + \tr h_\epsilon(t)^{-1}))\leq |S_\epsilon(t)|_{\hat{H}}\leq r^{\frac{1}{2}}\log (\tr h_\epsilon(t) + \tr h_\epsilon(t)^{-1}),
\end{equation}
where $r=\rank (E)$.
By (\ref{H001}) and (\ref{220}), we have
\begin{equation}\label{C0a}
\int_{\tilde{M}}\log(\tr h_\epsilon(t)+\tr h_\epsilon^{-1}(t))-\log(2r)\frac{\omega_{\epsilon}^{n}}{n!}\leq \hat{C}_{1}t,
\end{equation}
and
\begin{equation}\label{C0b}
\log(\tr h_\epsilon(t)+\tr h_\epsilon^{-1}(t))-\log(2r)\leq 2\tilde{C}_{1}(\delta^{-1})T
\end{equation}
for all $(x, t) \in (\tilde{M}\setminus B_{\omega_1}(\delta))\times [0, T]$. Then, we have the following local $C^{0}$-estimate of $h_{\epsilon}(t)$.

\medskip

\begin{lemma}\label{lem 2.3}
There exists a constant $\overline{C}_{0}(\delta^{-1}, T)$ which is independent of $\epsilon$ such that
\begin{equation}\label{C01}
|S_{\epsilon}(t)|_{\hat{H}}(x)\leq \overline{C}_{0}(\delta^{-1}, T)
\end{equation}
for all $(x, t) \in (\tilde{M}\setminus B_{\omega_1}(\delta))\times [0, T]$, and all $0< \epsilon \leq 1$.
\end{lemma}

\medskip

In the following lemma, we derive a local $C^{1}$-estimate of $h_{\epsilon }(t)$.

\medskip

\begin{lemma}\label{lem 2.4}
Let $T_\epsilon(t)=h_\epsilon^{-1}(t)\partial_{\hat{H}}h_\epsilon(t)$. Assume that there exists a constant
$\overline{C}_{0}$ such that
\begin{equation}\label{C02}
\max_{(x, t)\in (\tilde{M}\setminus B_{\omega_1}(\delta))\times [0, T]}|S_{\epsilon}(t)|_{\hat{H}}(x)\leq \overline{C}_{0},
\end{equation}
for all $0< \epsilon \leq 1$. Then, there exists a constant $\overline{C}_{1}$  depending only on $\overline{C}_{0}$ and $\delta ^{-1}$ such that
\begin{equation}\label{C11}
\max_{(x, t)\in (\tilde{M}\setminus B_{\omega_1}(\frac{3}{2}\delta))\times [0, T]}|T_{\epsilon}(t)|_{\hat{H}, \omega_{\epsilon}}\leq \overline{C}_{1}
\end{equation}
for all $0< \epsilon \leq 1$.
\end{lemma}

\medskip

{\bf Proof. } By a direct calculation, we have
\begin{equation}\label{c1002}
\begin{split}
&(\Delta_\epsilon-\frac{\partial}{\partial t})\tr h_\epsilon(t)\\
=&\ 2 \tr (-\sqrt{-1}\Lambda_{\omega_{\epsilon}}\overline{\partial}h_\epsilon(t)\cdot h_\epsilon^{-1}(t)\cdot\partial_{\hat{H}}h_\epsilon(t))+2\tr(h_{\epsilon}(t)\Phi (\hat{H}, \omega_{\epsilon})) \\
 & +2\sqrt{-1}\Lambda_{\omega_{\epsilon}}\tr \{h_{\epsilon}(t)\circ ([\phi , \phi ^{\ast H_{\epsilon}(t)}]-[\phi , \phi ^{\ast \hat{H}}])\}\\
=&\ 2 \tr (-\sqrt{-1}\Lambda_{\omega_{\epsilon}}\overline{\partial}h_\epsilon(t)\cdot h_\epsilon^{-1}(t)\cdot\partial_{\hat{H}}h_\epsilon(t))+2\tr(h_{\epsilon}(t)\Phi (\hat{H}, \omega_{\epsilon})) \\
 & +2\sqrt{-1}\Lambda_{\omega_{\epsilon}}\tr \{[\phi , h_{\epsilon}(t)]\wedge h_{\epsilon}^{-1}(t)[h_{\epsilon}(t) ,  \phi ^{\ast \hat{H}}]\}\\
\geq &\ 2 \tr (-\sqrt{-1}\Lambda_{\omega_{\epsilon}}\overline{\partial}h_\epsilon(t)\cdot h_\epsilon^{-1}(t)\cdot\partial_{\hat{H}}h_\epsilon(t))+2\tr(h_{\epsilon}(t)\Phi (\hat{H}, \omega_{\epsilon})), \\
\end{split}
\end{equation}
\begin{equation}
\frac{\partial }{\partial t}T_\epsilon(t) =\partial_{H_{\epsilon}(t)}(h_\epsilon^{-1}(t)\frac{\partial }{\partial t}h_\epsilon(t))=-2\partial_{H_{\epsilon}(t)}(\Phi (H_{\epsilon}(t), \omega_{\epsilon})),
\end{equation}
and
\begin{equation}\label{c1001}
\begin{split}
&(\Delta_\epsilon-\frac{\partial}{\partial t})|T_\epsilon(t)|^2_{H_\epsilon(t), \omega_{\epsilon}} \geq 2|\nabla_{H_\epsilon(t)}T_\epsilon(t)|^2_{H_\epsilon(t), \omega_{\epsilon}}\\
& -\check{C}_1(|\Lambda_{\omega_{\epsilon}}F_{H_{\epsilon}(t)}|_{H_\epsilon(t)}+|F_{\hat{H}}|_{H_\epsilon(t), \omega_{\epsilon}}+|\phi|_{H_\epsilon(t), \omega_{\epsilon}}^{2}
+|Ric(\omega_{\epsilon})|_{\omega_{\epsilon}})|T_\epsilon(t)|^2_{H_\epsilon(t), \omega_{\epsilon}}\\
& -\check{C}_2|\nabla_{\hat{H}}(\Lambda_{\omega_{\epsilon}}F_{\hat{H}})|_{H_\epsilon(t), \omega_{\epsilon}}|T_\epsilon(t)|_{H_\epsilon(t), \omega_{\epsilon} }-|\nabla_{\hat{H}}\phi |_{H_\epsilon(t), \omega_{\epsilon}}^{2},\\
\end{split}
\end{equation}
where constants $\check{C}_1, \check{C}_2$ depend only on the dimension $n$ and the rank $r$.

By the local $C^{0}$-assumption (\ref{C02}), the local estimate (\ref{220}) and the definition of $\omega_\epsilon$, it is easy to see that all coefficients in the right term of (\ref{c1001}) are uniformly local bounded outside $\tilde{\Sigma }$. Then there exists a constant $\check{C}_3$ depending only on $\delta^{-1}$ and $\overline{C}_{0}$ such that
\begin{equation}\label{c10003}
\begin{split}
(\Delta_\epsilon-\frac{\partial}{\partial t})|T_\epsilon(t)|^2_{H_\epsilon(t), \omega_{\epsilon}} \geq&\  2|\nabla_{H_\epsilon(t)}T_\epsilon(t)|^2_{H_\epsilon(t), \omega_{\epsilon}}\\
& -\check{C}_3|T_\epsilon(t)|^2_{H_\epsilon(t), \omega_{\epsilon}}-\check{C}_3 \\
\end{split}
\end{equation}
on the domain $\tilde{M}\setminus B_{\omega_1}(\delta)\times [0, T]$.

Let $\varphi_1$, $\varphi_2$ be nonnegative cut-off functions satisfying:
\begin{equation}
\varphi_1(x)=\left\{ \begin{array}{ll}
0, &   x\in B_{\omega_1}(\frac{5}{4}\delta),\\
1, &   x\in \tilde{M}\setminus B_{\omega_1}(\frac{3}{2}\delta),
\end{array} \right.
\end{equation}

\begin{equation}
\varphi_2(x)=\left\{ \begin{array}{ll}
0, &   x\in B_{\omega_1}(\delta),\\
1, &   x\in \tilde{M}\setminus B_{\omega_1}(\frac{5}{4}\delta),
\end{array} \right.
\end{equation}
and $|d\varphi_i|_{\omega_1}^2\leq \frac{8}{\delta^2}$, $-\frac{c}{\delta^2}\omega_1\leq \sqrt{-1}\partial\bar\partial\varphi_i\leq \frac{c}{\delta^2}\omega_1$.
By the inequality (\ref{metric02}), there exists a constant $C_1(\delta^{-1})$ depending only on $\delta^{-1}$ such that
\begin{equation}
(|d\varphi_i|^2_{\omega_\epsilon}+|\Delta_\epsilon\varphi_i|)\leq C_1(\delta^{-1}),
\end{equation}
for all $0< \epsilon \leq 1$.

We consider the following test function
\begin{equation}
f(\cdot, t)=\varphi_1^2|T_\epsilon(t)|_{H_\epsilon(t), \omega_{\epsilon}}^2+W\varphi_2^2 \tr h_\epsilon(t),
\end{equation}
where the constant $W$ will be chosen large enough later. From (\ref{c1002}) and (\ref{c1001}), we have
\begin{equation}
\begin{split}
&(\Delta_\epsilon-\frac{\partial}{\partial t})f\\
=&\ \varphi_1^2(2|\nabla_{H_\epsilon(t)}T_\epsilon(t)|^2_{H_\epsilon(t), \omega_{\epsilon}}
-\check{C}_3|T_\epsilon(t)|^2_{H_\epsilon(t), \omega_{\epsilon}}-\check{C}_3+\Delta_{\omega_{\epsilon}}\varphi_1^2|T_\epsilon(t)|^2_{H_\epsilon(t), \omega_{\epsilon}}\\
&+4\langle\varphi_1\nabla\varphi_1, \nabla|T_\epsilon(t)|_{H_\epsilon(t), \omega_{\epsilon}}^2\rangle_{\omega_{\epsilon}}+W\Delta_{\omega_{\epsilon}}\varphi_2^2\tr h_\epsilon(t)+ 4W\langle\varphi_2\nabla\varphi_2, \nabla\tr h_\epsilon(t)\rangle_{\omega_{\epsilon}}\\
&+2W\varphi_2^2(\tr (\sqrt{-1}\Lambda_{\omega_\epsilon}h_\epsilon^{-1}(t)\partial_{\hat{H}}h_\epsilon(t)\bar{\partial}h_\epsilon(t)))+\tr (h_\epsilon(t)(\Phi (\hat{H}, \omega_{\epsilon}))).\\
\end{split}
\end{equation}
We use
\begin{equation}
\begin{split}
2\langle \varphi_1\nabla\varphi_1, \nabla|T_\epsilon(t)|_{H_\epsilon(t), \omega_{\epsilon}}^2 \rangle_{\omega_\epsilon}&\geq  -4\varphi_1|\nabla\varphi_1|_{\omega_\epsilon}|T_\epsilon(t)|_{H_\epsilon(t), \omega_{\epsilon}}|\nabla_{H_\epsilon(t)}T_\epsilon(t)|_{H_\epsilon(t), \omega_{\epsilon}}\\
&\geq -\varphi_1^2|T_\epsilon(t)|_{H_\epsilon(t), \omega_{\epsilon}}^2-4|\nabla\varphi_1|_{\omega_{\epsilon}}^2|T_\epsilon(t)|_{H_\epsilon(t), \omega_{\epsilon}}^2,\\
\end{split}
\end{equation}
\begin{equation}
W\langle \varphi_2\nabla\varphi_2, \nabla\tr h_\epsilon(t) \rangle_{\omega_\epsilon} \geq -\varphi_2^2|\nabla\tr h_\epsilon(t)|_{H_\epsilon(t), \omega_{\epsilon}}^2-W^{2}|\nabla\varphi_2|_{\omega_\epsilon}^2,
\end{equation}
and
\begin{equation}
\begin{split}
&|T_\epsilon(t)|^2_{H_\epsilon(t), \omega_{\epsilon }}\\
=&\ \tr(\sqrt{-1}\Lambda_{\omega_\epsilon}h_\epsilon^{-1}(t)\partial_{\hat{H}}
h_\epsilon(t)H_\epsilon^{-1}(t)\overline{(h_\epsilon^{-1}(t)\partial_{\hat{H}}h_\epsilon(t))}^TH_\epsilon(t))\\
=&\ \tr(\sqrt{-1}\Lambda_{\omega_\epsilon}h_\epsilon^{-1}(t)\partial_{\hat{H}}h_\epsilon(t)h_\epsilon^{-1}(t)\bar{\partial}h_\epsilon(t))\\
\leq &\  e^{\overline{C}_{0}}\tr(\sqrt{-1}\Lambda_{\omega_\epsilon}h_\epsilon^{-1}(t)\partial_{\hat{H}}h_\epsilon(t)\bar{\partial}h_\epsilon(t)),\\
\end{split}
\end{equation}
and choose
\begin{equation}
W=(\check{C}_3+4C_1(\delta^{-1})+2r)e^{\overline{C}_{0}}+1.
\end{equation}
Then there exists a positive constant $\tilde{C}_0$ depending only on $\overline{C}_{0}$ and $\delta ^{-1}$ such that
\begin{equation}\label{c101}
(\Delta_\epsilon-\frac{\partial}{\partial t})f\geq \varphi_1^2|\nabla_{H_\epsilon(t)}T_\epsilon(t)|^2_{H_\epsilon(t), \omega_{\epsilon}}+\varphi_2^2|T_\epsilon(t)|^2_{H_\epsilon(t), \omega_{\epsilon}}-\tilde{C}_0
\end{equation}
on $\tilde{M}\times [0, T]$. Let $f (q, t_{0})=\max_{\tilde{M}\times [0, T]}\eta$, by the definition of $\varphi_{i}$ and the uniform local $C^{0}$-assumption of $h_{\epsilon}(t)$,  we can suppose that:
$$
(q, t_{0}) \in \tilde{M}\setminus B_{\omega_1}(\frac{5}{4}\delta)\times (0, T].
$$
By the inequality (\ref{c101}), we have
\begin{equation}
|T_\epsilon(t_{0})|^2_{H_\epsilon(t_{0}), \omega_{\epsilon }}(q)\leq \tilde{C}_0.
\end{equation}
So there exists a constant $\overline{C}_{1}$ depending only on $\overline{C}_{0}$ and $\delta ^{-1}$, such that
\begin{equation}
|T_\epsilon(t)|^2_{H_\epsilon(t), \omega_\epsilon}(x)\leq \overline{C}_{1}
\end{equation}
for all $(x, t)\in \tilde{M}\setminus B_{\omega_1}(\frac{3}{2}\delta)\times [0, T]$ and all $0< \epsilon \leq 1$.

\hfill $\Box$ \\

          One can get the local uniform $C^{\infty}$ estimates of $h_{\epsilon}(t)$ by the standard Schauder estimate of the parabolic equation  after getting the local $C^{0}$ and $C^1$ estimates. But by applying the parabolic Schauder estimates, one can only get the  uniform $C^{\infty}$ estimates of $h_{\epsilon} (t)$ on $\tilde{M}\setminus B_{\omega_1}(\delta)\times [\tau , T]$, where $\tau >0$ and the uniform estimates depend on $\tau ^{-1}$. In the following, we first use the maximum principle to get a local uniform bound on the curvature $|F_{H_{\epsilon}(t)}|_{H_{\epsilon}(t), \omega_{\epsilon }}$, then we apply the elliptic estimates to get local uniform $C^{\infty}$ estimates. The benefit of our argument is that we can get   uniform $C^{\infty}$ estimates of $h_{\epsilon}(t)$ on $\tilde{M}\setminus B_{\omega_1}(\delta)\times [0, T]$.  In the following, for simplicity, we denote
 \begin{equation}\Xi_{\epsilon, j}=|\nabla_{H_{\epsilon }(t)}^{j}(F_{H_{\epsilon}(t)}+[\phi ,
\phi^{\ast H_{\epsilon}(t)}])|_{H_{\epsilon}(t), \omega_{\epsilon}}^{2}(x)+|\nabla_{H_{\epsilon}(t)}^{j+1} \phi |_{H_{\epsilon}(t), \omega_{\epsilon}}^{2}\end{equation}
for $j=0, 1, \cdots $. Here $\nabla_{H_{\epsilon}(t)}$ denotes the covariant derivative with respect to the Chern connection $D_{H_{\epsilon}(t)}$ of $H_{\epsilon}(t)$ and the Riemannian connection $\nabla_{\omega_{\epsilon}}$ of $\omega_{\epsilon}$.

\medskip

\begin{lemma}\label{lem 2.5}
Assume that there exists a constant
$\overline{C}_{0}$ such that
\begin{equation}
\max_{(x, t)\in (\tilde{M}\setminus B_{\omega_1}(\delta))\times [0, T]}|S_{\epsilon}(t)|_{\hat{H}}(x)\leq \overline{C}_{0},
\end{equation}
for all $0< \epsilon \leq 1$. Then, for every integer $k\geq 0$, there exists a constant $\overline{C}_{k+2}$  depending only on $\overline{C}_{0}$, $\delta ^{-1}$ and $k$,  such that
\begin{equation}\label{C11}
\max_{(x, t)\in (\tilde{M}\setminus B_{\omega_1}(2\delta))\times [0, T]}\Xi_{\epsilon, k}\leq \overline{C}_{k+2}
\end{equation}
for all $0< \epsilon \leq 1$. Furthermore, there exist  constants $\hat{C}_{k+2}$  depending only on $\overline{C}_{0}$, $\delta ^{-1}$ and $k$,  such that
\begin{equation}
\max_{(x, t)\in (\tilde{M}\setminus B_{\omega_1}(2\delta))\times [0, T]}|\nabla_{\hat{H}}^{k+2}h_{\epsilon}|_{\hat{H}, \omega_\epsilon}\leq \hat{C}_{k+2}
\end{equation}
for all $0< \epsilon \leq 1$.
\end{lemma}

\medskip

{\bf Proof. }
By computing, we have the following inequalities (see Lemma 2.4 and Lemma 2.5 in (\cite{LZ1}) for details):
\begin{equation}\label{cc01}
\begin{split}
&(\Delta_{\epsilon } -\frac{\partial }{\partial t})|\nabla _{H_{\epsilon}(t)}\phi
|_{H_{\epsilon}(t), \omega_\epsilon}^{2}-2|\nabla_{H_{\epsilon}(t)}\nabla_{H_{\epsilon}(t)}\phi |_{H_{\epsilon}(t), \omega_\epsilon}^{2}\\
\geq &-C_{7}(|F_{H_{\epsilon}(t)}|_{H_{\epsilon}(t), \omega_\epsilon}+|Rm(\omega_{\epsilon})|_{\omega_\epsilon}+|\phi |_{H_{\epsilon}(t), \omega_\epsilon}^{2})|\nabla_{H_{\epsilon}(t)} \phi
|_{H_{\epsilon}(t), \omega_\epsilon}^{2}\\
&-C_{7}|\phi |_{H_{\epsilon}(t), \omega_\epsilon}|\nabla Ric (\omega_{\epsilon})|_{\omega_\epsilon}|\nabla_{H_{\epsilon}(t)} \phi|_{H_{\epsilon}(t), \omega_\epsilon},\\
\end{split}
\end{equation}
\begin{equation}
\begin{split}
&(\Delta_{\epsilon} -\frac{\partial }{\partial t})|F_{H_{\epsilon}(t)}+[\phi ,
\phi^{\ast H_{\epsilon}(t)}]|_{H_{\epsilon}(t), \omega_\epsilon}^{2}-2|\nabla_{H_{\epsilon}(t)}(F_{H_{\epsilon}(t)}+[\phi , \phi^{\ast H_{\epsilon}(t)
}])|_{H_{\epsilon}(t), \omega_\epsilon}^{2} \\
 \geq &-C_{8}(|F_{H_{\epsilon}(t)}+[\phi , \phi^{\ast H_{\epsilon}(t)}]|_{H_{\epsilon}(t), \omega_\epsilon}^{2}+|\nabla_{H_{\epsilon}(t)} \phi |_{H_{\epsilon}(t), \omega_\epsilon}^{2})^{\frac{3}{2}}\\
&-C_{8} (|\phi
|_{H_{\epsilon}(t), \omega_\epsilon}^{2}+|Rm(\omega_{\epsilon})|_{\omega_\epsilon})(|F_{H_{\epsilon}(t)}+[\phi , \phi^{\ast H_{\epsilon}(t)}]|_{H_{\epsilon}(t), \omega_\epsilon}^{2}+|\nabla_{H_{\epsilon}(t)} \phi |_{H_{\epsilon}(t), \omega_\epsilon}^{2}),\\
\end{split}
\end{equation}
then
\begin{equation}\label{cc02}
\begin{split}
(\Delta_{\epsilon} -\frac{\partial }{\partial t})\Xi_{\epsilon, 0} \geq&\ 2\Xi_{\epsilon, 1}
-C_{8}(\Xi_{\epsilon, 0})^{\frac{3}{2}}\\
&-C_{8} (|\phi
|_{H_{\epsilon}(t), \omega_\epsilon}^{2}+|Rm(\omega_{\epsilon})|_{\omega_\epsilon})(\Xi_{\epsilon, 0})-C_{8}|\nabla Ric (\omega_{\epsilon })|^{2}_{\omega_{\epsilon}},\\
\end{split}
\end{equation}
where $C_{7}$, $C_{8}$ are constants depending only on the complex
dimension $n$ and the rank $r$.  Furthermore, we have
\begin{equation}\label{Fk}
\begin{split}
&(\Delta_{\epsilon}-\frac{\partial }{\partial t} )\Xi_{\epsilon, j}\\
\geq&\ 2\Xi_{\epsilon, j+1}-\acute{C}_{j}(\Xi_{\epsilon, j})^{\frac{1}{2}}\{\sum_{i+k=j}((\Xi_{\epsilon, i})^{\frac{1}{2}}+|\phi
|_{H_{\epsilon}(t), \omega_\epsilon}^{2}+|Rm(\omega_{\epsilon})|_{\omega_{\epsilon}}+|\nabla Ric (\omega_{\epsilon })|_{\omega_{\epsilon}})\\ & \cdot
((\Xi_{\epsilon, k})^{\frac{1}{2}}+|\phi
|_{H_{\epsilon}(t), \omega_\epsilon}^{2}+|Rm(\omega_{\epsilon})|_{\omega_{\epsilon}}+|\nabla Ric (\omega_{\epsilon })|_{\omega_{\epsilon}})\},\\
\end{split}
\end{equation}
where $\acute{C}_{j}$ is a positive constant depending only on the complex
dimension $n$, the rank $r$ and $j$.
Direct computations yield the following inequality (see (2.5) in (\cite{LZ1}) for details):
\begin{equation}\label{phi002}
\begin{split}
&(\Delta_{\epsilon}-\frac{\partial }{\partial t})|\phi|_{H_{\epsilon}(t),\omega_{\epsilon}}^{2}\geq 2|\nabla_{H_{\epsilon}(t)} \phi |_{H_{\epsilon}(t), \omega_{\epsilon}}^{2}\\
&+2|\Lambda_{\omega_{\epsilon}}[\phi , \phi^{\ast H_{\epsilon}(t)}]|_{H_{\epsilon}(t)}^{2}-2|Ric(\omega_{\epsilon})|_{\omega_{\epsilon}}|\phi|_{H_{\epsilon}(t),\omega_{\epsilon}}^{2}.\\
\end{split}
\end{equation}
From the local $C^{0}$-assumption (\ref{C02}), we see that $|\phi|_{H_{\epsilon}(t), \omega_{\epsilon}}$ is also uniformly bounded on $\tilde{M}\setminus B_{\omega_1}(\delta)\times [0, T]$. By Lemma \ref{lem 2.4}, we have $|T_{\epsilon }(t)|_{H_{\epsilon}(t), \omega_{\epsilon}}$ is uniformly bounded on $\tilde{M}\setminus B_{\omega_1}(\frac{3}{2}\delta)\times [0, T]$. We choose a constant $\hat{C}$ depending only on $\delta^{-1}$ and $\overline{C}_{0}$ such that
\begin{equation}
\frac{1}{2}\hat{C}\leq \hat{C} -(|\phi|_{H_{\epsilon}(t), \omega_{\epsilon}}^{2}+|T_{\epsilon }(t)|_{H_{\epsilon}(t), \omega_{\epsilon}}^{2})(x)\leq \hat{C}
\end{equation}
on $\tilde{M}\setminus B_{\omega_1}(\frac{3}{2}\delta)\times [0, T]$.
We consider the test function:
\begin{equation}
\zeta (x, t)=\rho^{2}\frac{\Xi_{\epsilon , 0}(x, t)}{\hat{C}-(|\phi|_{H_{\epsilon}(t), \omega_{\epsilon}}^{2}+|T_{\epsilon }(t)|_{H_{\epsilon}(t), \omega_{\epsilon}}^{2})(x)},
\end{equation}
where $\rho $ is a cut-off function satisfying:
\begin{equation}
\rho (x)=\left\{ \begin{array}{ll}
0, &   x\in B_{\omega_1}(\frac{13}{8}\delta),\\
1, &   x\in \tilde{M}\setminus B_{\omega_1}(\frac{7}{4}\delta),
\end{array} \right.
\end{equation}
and $|d\rho |_{\omega_1}^2\leq \frac{8}{\delta^2}$, $-\frac{c}{\delta^2}\omega_1\leq \sqrt{-1}\partial\bar\partial\rho \leq \frac{c}{\delta^2}\omega_1$. We suppose $(x_{0}, t_{0})\in \tilde{M}\setminus B_{\omega_1}(\frac{3}{2}\delta)\times (0, T]$ is a maximum point of $\zeta $. Using (\ref{c10003}), (\ref{cc01}), (\ref{cc02}), (\ref{phi002}) and the fact $\nabla\zeta =0$ at the point $(x_{0}, t_{0})$, we have
\begin{equation}
\begin{split}
0 \geq&\ (\Delta_{\epsilon }-\frac{\partial }{\partial t})\zeta |_{(x_{0}, t_{0})}\\
=&\ \frac{1}{\hat{C}-(|\phi|_{H_{\epsilon}(t), \omega_\epsilon}^{2}+|T_{\epsilon }(t)|_{H_{\epsilon}(t), \omega_\epsilon}^{2})} (\Delta_{\epsilon }-\frac{\partial }{\partial t})(\rho^{2}\Xi_{\epsilon , 0}) \\& -\rho^{2}\frac{\Xi_{\epsilon , 0}}{(\hat{C}-(|\phi|_{H_{\epsilon}(t), \omega_\epsilon}^{2}+|T_{\epsilon }(t)|_{H_{\epsilon}(t), \omega_\epsilon}^{2}))^{2}}(\Delta_{\epsilon }-\frac{\partial }{\partial t})(\hat{C}-(|\phi|_{H_{\epsilon}(t), \omega_\epsilon}^{2}+|T_{\epsilon }(t)|_{H_{\epsilon}(t), \omega_\epsilon}^{2}))\\
& -\frac{2}{\hat{C}-(|\phi|_{H_{\epsilon}(t), \omega_\epsilon}^{2}+|T_{\epsilon }(t)|_{H_{\epsilon}(t), \omega_\epsilon}^{2})}\nabla (\zeta)\cdot  \nabla (\hat{C}-(|\phi|_{H_{\epsilon}(t), \omega_\epsilon}^{2}+|T_{\epsilon }(t)|_{H_{\epsilon}(t), \omega_\epsilon}^{2}))\\
\geq&\ \frac{\Xi_{\epsilon , 0}}{(\hat{C}-(|\phi|_{H_{\epsilon}(t), \omega_\epsilon}^{2}+|T_{\epsilon }(t)|_{H_{\epsilon}(t), \omega_\epsilon}^{2}))^{2}}\{\rho^{2}\frac{2\Xi_{\epsilon , 0}-\check{C}_3 |T_\epsilon(t)|^2_{H_\epsilon(t), \omega_\epsilon}-\check{C}_3}{\hat{C}-(|\phi|_{H_{\epsilon}(t), \omega_\epsilon}^{2}+|T_{\epsilon }(t)|_{H_{\epsilon}(t), \omega_\epsilon}^{2})}\\
& -\rho^{2}\frac{2|Ric(\omega_{\epsilon})|_{\omega_\epsilon}|\phi|_{H_{\epsilon}(t), \omega_\epsilon}^{2}}{\hat{C}-(|\phi|_{H_{\epsilon}(t), \omega_\epsilon}^{2}+|T_{\epsilon }(t)|_{H_{\epsilon}(t), \omega_\epsilon}^{2})}\\
& -C_{8}\rho^{2} \Xi_{\epsilon , 0}^{\frac{1}{2}}-C_{8}\rho^{2} (|\phi
|_{H_{\epsilon}(t), \omega_\epsilon}^{2}+|Rm(\omega_{\epsilon})|_{\omega_\epsilon})-8|d\rho |_{\omega_{\epsilon}}^{2}+\Delta_{\omega_{\epsilon}}\rho^{2}\}\\
& -C_{8}\frac{\rho^{2}|\nabla Ric (\omega_{\epsilon})|_{\omega_{\epsilon}}^{2}}{\hat{C}-(|\phi|_{H_{\epsilon}(t), \omega_\epsilon}^{2}+|T_{\epsilon }(t)|_{H_{\epsilon}(t), \omega_\epsilon}^{2})}.\\
\end{split}
\end{equation}
So there exist positive constants $\dot{C}_{2}$ and $\overline{C}_{2}$ depending only on $\overline{C}_{0}$ and $\delta ^{-1}$, such that
\begin{equation}
\zeta (x_{0}, t_{0})\leq \dot{C}_{2},
\end{equation}
and
\begin{equation}
\Xi_{\epsilon , 0} (x, t)\leq \overline{C}_{2}
\end{equation}
for all $(x, t)\in \tilde{M}\setminus B_{\omega_1}(\frac{7}{4}\delta)\times [0, T]$.

Furthermore, we choose two suitable cut-off functions $\rho_{1}$, $\rho_{2}$, a suitable constant $A$ which depends only on $\overline{C}_{0}$ and $\delta ^{-1}$, and  a test function
\begin{equation}
\zeta_{1}(x, t)=\rho_{1}^{2}\Xi_{\epsilon , 1 }+ A \rho_{2}^{2}\Xi_{\epsilon , 0}.
\end{equation}
Running a similar argument as above, we can show that there exist constants $\overline{C}_{3}$ and $\dot{C}_{3}$ depending only on $\overline{C}_{0}$ and $\delta ^{-1}$ such that
\begin{equation}
\Xi_{\epsilon , 1} (x, t)\leq \overline{C}_{3},
\end{equation}
and
\begin{equation}\label{c201}
|\nabla_{\hat{H}}F_{H_{\epsilon}(t)}|_{\hat{H}, \omega_{\epsilon}}^{2}\leq \dot{C}_{3}
\end{equation}
for all $(x, t)\in \tilde{M}\setminus B_{\omega_1}(\frac{15}{8}\delta)\times [0, T]$.

\medskip

Recalling the equality
 \begin{equation}
 \overline{\partial }\partial _{\hat{H}}h_{\epsilon}(t)=h_{\epsilon}(t)(F_{H_{\epsilon}(t)}-F_{\hat{H}})+\overline{\partial }h_{\epsilon}(t)\wedge (h_{\epsilon}(t))^{-1}\partial _{\hat{H}}h_{\epsilon}(t)
 \end{equation}
 and noting that K\"ahler metrics $\omega_{\epsilon}$ are uniform locally quasi-isometry to $\pi^{\ast }\omega $ outside the exceptional divisor $\tilde{\Sigma}$, by standard elliptic estimates, because we have local uniform bounds on $h_{\epsilon }$, $T_{\epsilon}$, $F_{H_{\epsilon}}$ and $ F_{\hat{H}}$,  we get a uniform $C^{1, \alpha }$-estimate of $h_{\epsilon }$ on $\tilde{M}\setminus B_{\omega_1}(\frac{61}{32}\delta)\times [0, T]$.
 
 \medskip

 We can iterate this procedure by induction and then obtain  local uniform bounds for $\Xi_{\epsilon , k}$, $|\nabla_{\hat{H}}^{k}F_{H_{\epsilon}(t)}|_{\hat{H}, \omega_\epsilon}^{2}$, and $\|h_{\epsilon}\|_{C^{k+1,\alpha}}$ on $\tilde{M}\setminus B_{\omega_1}(2\delta)\times [0, T]$ for any $k\geq 1$.

\hfill $\Box$ \\

\medskip

From the above local uniform $C^{\infty}$-bounds on $H_{\epsilon}$, we get the following Lemma.

\medskip

\begin{lemma}\label{lem 2.6}
By choosing a subsequence,   $H_{\epsilon}( t)$  converges to $H(x, t)$ locally in $C^\infty$ topological on $\tilde{M}\setminus \tilde{\Sigma }\times [0, \infty)$ as $\epsilon\rightarrow 0$ and $H(t)$ satisfies (\ref{SSS1}).
\end{lemma}
\medskip

\section{Uniform estimate of the Higgs field }
\setcounter{equation}{0}

In this section, we prove that the norm $|\phi |_{H(t), \omega}$ is uniformly bounded along the heat flow (\ref{SSS1}) for $t\geq t_{0}>0$.

Firstly, we know $|\phi|_{\hat{H}, \omega_\epsilon}^2
\in L^{1}(\tilde{M}, \omega_\epsilon )$ and the $L^{1}$-norm is uniformly bounded.
In fact,
\begin{eqnarray}
\begin{array}{lll}
\int_{\tilde{M}}|\phi|_{\hat{H}, \omega_\epsilon}^2\frac{\omega_\epsilon^n}{n!}&=\int_{\tilde{M}}\tr(\sqrt{-1}\Lambda_{\omega_\epsilon}(\phi\wedge\phi^{\ast\hat{H}}))\frac{\omega_\epsilon^n}{n!}\\
&=\int_{\tilde{M}}\tr (\phi\wedge\phi^{\ast\hat{H}})\wedge\frac{\omega_\epsilon^{n-1}}{(n-1)!}\leq \check{C}_{\phi}< \infty,
\end{array}
\end{eqnarray}
where $\check{C}_{\phi}$ is a positive constant independent of $\epsilon$. Moreover, we will show the $L^{1+2a}$-norm of $|\phi|_{\hat{H}, \omega_\epsilon}^2$ is also uniformly bounded, for any $0\leq 2a< \frac{1}{2}$.
Let's recall Lemma 5.5 in \cite{Sib} (see also Lemma 5.8 in \cite{LZ3}).

\medskip

\begin{lemma}\label{lem 3.1}{\bf(\cite{Sib})}
Let $(M, \omega)$ be a compact K\"ahler manifold of complex dimension $n$, and $\pi: \tilde{M}\rightarrow M$ be a blow-up along a smooth complex sub-manifold $\Sigma$ of complex codimension $k$ where $k\geq 2$. Let $\eta$ be a K\"ahler metric on $\tilde{M}$, and consider the family of K\"ahler metric $\omega_\epsilon=\pi^\ast\omega+\epsilon\eta$. Then for any $0\leq 2a <\frac{1}{k-1}$, we have $\frac{\eta^n}{\omega_\epsilon^n}\in L^{2a}(\tilde{M}, \eta)$, and the $L^{2a}(\tilde{M}, \eta)$-norm of $\frac{\eta^n}{\omega_\epsilon^n}$ is uniformly bounded independent of $\epsilon$, i.e. there is a positive constant $C^\ast$ such that
\begin{equation}\label{Laa01}
\int_{\tilde{M}}(\frac{\eta^n}{\omega_\epsilon^n})^{2a}\frac{\eta^n}{n!}\leq C^\ast
\end{equation}
for all $0<\epsilon \leq 1$.
\end{lemma}

\medskip

Since $\phi \in \Omega^{1,0 }(\End (E))$ is a smooth section and $\omega_\epsilon=\pi^\ast\omega+\epsilon\eta$, there exists a uniform constant $\tilde{C}_{\phi }$ such that
\begin{equation}
\Big(\frac{|\phi|_{\hat{H}, \omega_\epsilon}^2\frac{\omega_\epsilon^n}{n!}}{\frac{\eta^n}{n!}}\Big)=\frac{n\tr (\phi\wedge\phi^{\ast\hat{H}})\wedge\omega_\epsilon^{n-1}}{\eta^n}\leq \tilde{C}_{\phi }
\end{equation}
for all $0<\epsilon \leq 1$.
 By (\ref{Laa01}), for any $0\leq 2a <\frac{1}{2}$, there exists a uniform constant $C_{\phi }$ such that
\begin{equation}\label{3.9}
\begin{split}
&\int_{\tilde{M}}|\phi|_{\hat{H}, \omega_\epsilon}^{2(1+2a)}\frac{\omega_\epsilon^n}{n!}\\
=&\int_{\tilde{M}}\Big(\frac{|\phi|_{\hat{H}, \omega_\epsilon}^2\frac{\omega_\epsilon^n}{n!}}{\frac{\eta^n}{n!}}\Big)^{1+2a}\Big(\frac{\eta^n}{\omega_\epsilon^n}\Big)^{1+2a}\frac{\omega_\epsilon^n}{n!}\\
=&\int_{\tilde{M}}\Big(\frac{|\phi|_{\hat{H}, \omega_\epsilon}^2\frac{\omega_\epsilon^n}{n!}}{\frac{\eta^n}{n!}}\Big)^{1+2a}\Big(\frac{\eta^n}{\omega_\epsilon^n}\Big)^{2a}\frac{\eta^n}{n!}\\
\leq &C_\phi
\end{split}
\end{equation}
for all $0<\epsilon \leq 1$. By limiting (\ref{3.9}), we have the following lemma.

\medskip

\begin{lemma}\label{lem 3.2}For any $0\leq 2a <\frac{1}{2}$, we have $|\phi|_{\hat{H}, \omega }^{2}\in L^{1+2a}(M\setminus \Sigma , \omega )$, i.e. there exists a  constant $C_{\phi }$ such that
\begin{equation}\label{Lp01}
\int_{M\setminus \Sigma }|\phi|_{\hat{H}, \omega }^{2(1+2a)}\frac{\omega ^n}{n!}\leq C_\phi .
\end{equation}
\end{lemma}

\medskip

On $M\setminus\Sigma$, we get ((2.5) in \cite{LZ1} for details)
\begin{equation}\label{eqn:1}
(\Delta-\frac{\partial}{\partial t})|\phi|_{H(t), \omega}^2\geq 2|\nabla_{H(t)}\phi|_{H(t), \omega}^2+2|\sqrt{-1}\Lambda_\omega[\phi, \phi^{\ast H(t)}]|_{H(t)}^2-2| Ric_\omega |_{\omega} |\phi|_{H(t), \omega}^2.
\end{equation}
By a direct computation, we have
\begin{equation}
\begin{aligned}
(\Delta-\frac{\partial}{\partial t})\log(|\phi|_{H(t), \omega}^2+e)=&\ \frac{1}{\log(|\phi|_{H(t), \omega}^2+e)}(\Delta-\frac{\partial}{\partial t})|\phi|_{H(t), \omega}^2-\frac{\nabla|\phi|_{H(t), \omega}^2\cdot\nabla|\phi|_{H(t), \omega}^2}{(|\phi|_{H(t), \omega}^2+e)^2}\\
\geq&\ \frac{1}{\log(|\phi|_{H(t), \omega}^2+e)}(\Delta-\frac{\partial}{\partial t})|\phi|_{H(t), \omega}^2-\frac{2|\nabla_{H(t)}^{1,0}\phi|_{H(t), \omega}^2\cdot|\phi|_{H(t), \omega}^2}{(|\phi|_{H(t), \omega}^2+e)^2}.
\end{aligned}
\end{equation}
Combining this with \eqref{eqn:1}, we obtain
\begin{equation}\label{eqn:3}
(\Delta-\frac{\partial}{\partial t})\log(|\phi|_{H(t), \omega}^2+e)\geq\frac{2|\Lambda_\omega[\phi, \phi^{\ast H(t)}]|_{H(t)}^2}{|\phi|_{H(t), \omega}^2+e}-2|Ric_\omega|_{\omega}
\end{equation}
on $M\setminus\Sigma$. Based on Lemma 2.7 in \cite{Si2}, we obtain
\begin{equation}\label{eqn:2}
|\sqrt{-1}\Lambda_\omega[\phi, \phi^{\ast H(t)}]|_{H(t)}=|[\phi, \phi^{\ast H(t)}]|_{H(t), \omega}\geq a_1|\phi|_{H(t), \omega}^2-a_2(|\phi|_{\hat{H}, \omega}^2+1),
\end{equation}
where $a_1$ and $a_2$ are positive constants depending only on $r$ and $n$. Then, for any $0\leq 2a < \frac{1}{2}$, we have
\begin{equation}
\begin{split}
 & 2|\Lambda_\omega[\phi, \phi^{\ast H(t)}]|_{H(t)}^2\\
\geq\ &(|\Lambda_\omega[\phi, \phi^{\ast H(t)}]|_{H(t)}+e)^2-6e^2\\
\geq\ &(|\Lambda_\omega[\phi, \phi^{\ast H(t)}]|_{H(t)}+e)^{1+\frac{a}{2}}-6e^2\\
\geq\ & a_{3}(|\phi|_{H(t), \omega }^2+e)^{1+\frac{a}{2}}-a_{4}|\phi|_{\hat{H}, \omega}^{2+a}-a_{5},
\end{split}
\end{equation}
where $a_{3}$, $a_{4}$ and $a_{5}$ are positive constants depending only on $a$, $r$ and $n$.
Then it is clear that (\ref{eqn:3}) implies:
\begin{equation}
(\Delta-\frac{\partial}{\partial t})\log(|\phi|_{H(t), \omega }^2+e)\geq a_3(|\phi|_{H(t), \omega}^2+e)^{\frac{a}{2}}-a_4|\phi|_{\hat{H}, \omega}^{2+a}-a_5-2|Ric_\omega|_{\omega},
\end{equation}
on $M\setminus\Sigma$.

\medskip

In the following, we denote:
\begin{equation}
f=\log(|\phi|_{H(t), \omega}^2+e).
\end{equation}
For any $b> 1$, we have:
\begin{equation}
\begin{aligned}
(\Delta-\frac{\partial}{\partial t})f^b=&\ bf^{b-1}(\Delta-\frac{\partial}{\partial t})f+b(b-1)|\nabla f|_{\omega}^2f^{b-2}\\
\geq&\ a_3bf^{b-1}(|\phi|_{H(t), \omega}^2+e)^{\frac{a}{2}}-a_4bf^{b-1}|\phi|_{\hat{H}, \omega}^{2+a}-(a_5+2|Ric_\omega|_{\omega})bf^{b-1}\\
& +b(b-1)|\nabla f|_{\omega}^2f^{b-2}.
\end{aligned}
\end{equation}

Choosing a  cut-off function $\varphi_{\delta }$ with
\begin{equation}
\varphi_\delta(x)=\left\{ \begin{array}{ll}
1, &   x\in M\setminus B_{2\delta}(\Sigma),\\
0, &   x\in B_{\delta}(\Sigma),
\end{array} \right.
\end{equation}
where $B_{\delta}=\{x\in M|d_{\omega}(x, \Sigma )<\delta \}$, and integrating by parts, we have
\begin{equation}
\begin{split}
 &-\frac{\partial}{\partial t}\int_M\varphi_\delta^4f^b\frac{\omega^n}{n!}=\int_M\varphi_\delta^4(\Delta-\frac{\partial}{\partial t})f^b\frac{\omega^n}{n!}+\int_M4\varphi_\delta^3\nabla\varphi_\delta\nabla f^b\frac{\omega^n}{n!}\\
\geq &\int_Ma_3b\varphi_\delta^4f^{b-1}(|\phi|_{H(t), \omega }^2+e)^{\frac{a}{2}}\frac{\omega^n}{n!}-\int_Ma_4b\varphi_\delta^4f^{b-1}
|\phi|_{\hat{H}, \omega}^{2+a}\frac{\omega^n}{n!}\\
 & -\int_M(a_5+2|Ric_\omega|_{\omega})b\varphi_\delta^4f^{b-1}\frac{\omega^n}{n!}+\int_Mb(b-1)\varphi_\delta^4
|\nabla f|_{\omega}^2f^{b-2}\frac{\omega^n}{n!}\\
 & -\int_M4b\varphi_\delta^3|\nabla\varphi_\delta|_{\omega}\cdot|\nabla f|_{\omega}f^{b-1}\frac{\omega^n}{n!}\\
 \geq&\int_Ma_3b\varphi_\delta^4f^{b-1}(|\phi|_{H(t), \omega }^2+e)^{\frac{a}{2}}\frac{\omega^n}{n!}-\int_Ma_4b\varphi_\delta^4f^{b-1}
(|\phi|_{\hat{H}, \omega}^{2})^{1+\frac{a}{2}}\frac{\omega^n}{n!}\\
 & -\int_M(a_5+2|Ric_\omega|_{\omega})b\varphi_\delta^4f^{b-1}\frac{\omega^n}{n!}
-\int_M\frac{4b}{b-1}\varphi_\delta^2|\nabla\varphi_\delta|_{\omega}^2f^b\frac{\omega^n}{n!}\\
 \geq &\int_Ma_3b\varphi_\delta^4f^{b-1}f^{(b-1)B}\frac{(|\phi|_{H(t), \omega}^2+e)^{\frac{a}{2}}}{f^{(b-1)B}}\frac{\omega^n}{n!}\\
 & -a_4b\Big(\int_M(\varphi_\delta^3f^{b-1})^p\frac{\omega^n}{n!}\Big)^{\frac{1}{p}}\Big(\int_M\varphi_\delta^q(|\phi|_{\hat{H}, \omega }^2)^{1+2a}
\frac{\omega^n}{n!}\Big)^{\frac{1}{q}}\\
 & -\int_M(a_5+2|Ric_\omega|_{\omega})b\varphi_\delta^4f^{b-1}\frac{\omega^n}{n!}\\
 & -\frac{4b}{b-1}\Big(\int_M\varphi_\delta^4f^{2b}\frac{\omega^n}{n!}\Big)^\frac{1}{2}
\Big(\int_M|\nabla\varphi_\delta|_{\omega}^4\frac{\omega^n}{n!}\Big)^{\frac{1}{2}},\\
\end{split}
\end{equation}
where $q=\frac{2(1+2a)}{2+a}$, $p=\frac{2(1+2a)}{3a}$ and $B=\frac{2(1+2a)}{3a}+\frac{2b}{b-1}$. We can see that there exists a constant $C(a,b)$ depending only on $a$ and $b$ such that
\begin{equation}
\frac{(|\phi|_{H(t), \omega}^2+e)^{\frac{a}{2}}}{(\log(|\phi|_{H(t), \omega}^2+e))^{(b-1)B}}\geq C(a,b).
\end{equation}
Since the complex codimension of $\Sigma$ is at least $3$, we can choose the cut-off function $\varphi_\delta$ such that
\begin{equation}\label{cut}
\int_M|\nabla\varphi_\delta|_{\omega}^4 \frac{\omega^n}{n!} \sim O(\delta^{-4}\delta^6)=O(\delta^2).
\end{equation}
By (\ref{Lp01}), we obtain
\begin{equation}\label{Lp04}
\begin{split}
-\frac{\partial}{\partial t}\int_M\varphi_\delta^4f^b\frac{\omega^n}{n!}\geq&  \ a_6\int_M\varphi_\delta^4f^{(b-1)B}\frac{\omega^n}{n!}-a_7\Big(\int_M\varphi_\delta^4f^{(b-1)B}\frac{\omega^n}{n!}\Big)^{\frac{1}{B}}\\
&-a_8\Big(\int_M\varphi_\delta^4f^{(b-1)B}\frac{\omega^n}{n!}\Big)^{\frac{1}{B}}
-a_{9}\Big(\int_M\varphi_\delta^4f^{(b-1)B}\frac{\omega^n}{n!}\Big)^{\frac{b}{(b-1)B}},\\
\end{split}
\end{equation}
where $a_i$ are positive constants depending only on $r, n, a, b, |Ric_\omega|_{\omega }, \Vol(M, \omega )$ and $C_{\phi}$ for $i=6, 7, 8, 9.$

\medskip

\begin{lemma}\label{lem 3.3}
For any $b>1$, there exists a constant $\hat{C}_{b}$ depending only on $r, n, b, |Ric_\omega|_{\omega}, \Vol(M, \omega )$ and $C_{\phi}$ such that
\begin{equation}\label{Lp02}
\int_{M\setminus \Sigma }(\log(|\phi|_{H(t), \omega }^2+e))^{b}\frac{\omega^n}{n!}\leq \hat{C}_{b}
\end{equation}
for all $t\geq 0$.
\end{lemma}

\medskip

{\bf Proof. } Suppose that $\int_M\varphi_\delta^4f^b\frac{\omega^n}{n!}(t^{\ast})=\max_{t\in [0, T]}\int_M\varphi_\delta^4f^b\frac{\omega^n}{n!}(t)$ with $t^{\ast }>0$. Choosing $a=\frac{1}{8}$ in (\ref{Lp04}), at point $t^{\ast}$, we have
\begin{equation}\label{Lp04}
\begin{split}
0\geq &-\frac{\partial}{\partial t}|_{t=t^{\ast}}\int_M\varphi_\delta^4f^b\frac{\omega^n}{n!}\\
\geq&\   a_6\int_M\varphi_\delta^4f^{(b-1)B}\frac{\omega^n}{n!}-a_7\Big(\int_M\varphi_\delta^4f^{(b-1)B}\frac{\omega^n}{n!}\Big)^{\frac{1}{B}}\\
&-a_8\Big(\int_M\varphi_\delta^4f^{(b-1)B}\frac{\omega^n}{n!}\Big)^{\frac{1}{B}}
-a_{9}\Big(\int_M\varphi_\delta^4f^{(b-1)B}\frac{\omega^n}{n!}\Big)^{\frac{b}{(b-1)B}}.\\
\end{split}
\end{equation}
This inequality implies that there exists a constant $\tilde{C}_{b}$ depending only on $r, n, b, |Ric_\omega|_{\omega}, \Vol(M, \omega )$ and $C_{\phi}$ such that
\begin{equation}
\int_M\varphi_\delta^4f^{(b-1)B}\frac{\omega^n}{n!}(t^{\ast})\leq \tilde{C}_{b}.
\end{equation}
So we have
\begin{equation}
\max_{t\in [0, T]}\int_M\varphi_\delta^4f^b\frac{\omega^n}{n!}(t)\leq \tilde{C}_{b}+\int_M(\log(|\phi|_{\hat{H}, \omega}^2+e))^b\frac{\omega^n}{n!}.
\end{equation}
Noting that the last term in the above inequality is also bounded, and letting $\delta \rightarrow 0$, we obtain the estimate (\ref{Lp02}) .

\hfill $\Box$ \\

By the heat equation (\ref{SSS1}), we have
\begin{equation}
|\frac{\partial}{\partial t}\log(|\phi|_{H(t), \omega}^2+e)|=\Big|\frac{\frac{\partial}{\partial t}|\phi|_{H(t), \omega}^2}{|\phi|_{H(t), \omega}^2+e}\Big|=\Big|\frac{2\langle[\Phi (H(t), \omega), \phi], \phi\rangle_{H(t)}}{|\phi|_{H(t), \omega}^2+e}\Big|\leq 2|\Phi (H(t), \omega)|_{H(t)},
\end{equation}
then
\begin{equation}
\Delta(\log(|\phi|_{H(t), \omega}^2+e))\geq -2|Ric_\omega|_{\omega}-2|\Phi (H(t), \omega)|_{H(t)}.
\end{equation}
 By (\ref{H005}), we have
 \begin{equation}\label{H00005}
 \max_{x\in M\setminus \Sigma }|\Phi (H(t), \omega) |_{H(t)}(x) \leq C_{K}(\tau )\hat{C}_{1} (t^{-n}+1).
\end{equation}
So there exists a positive constant $C^{\ast}(t_{0}^{-1})$ depending only on $t_{0}^{-1}$ and $|Ric_\omega|_{\omega}$ such that
\begin{equation}\label{e01}
\Delta(\log(|\phi|_{H(t), \omega}^2+e))\geq -C^{\ast}(t_{0}^{-1})
\end{equation}
on $M\setminus \Sigma $, for $t\geq t_{0}>0$.
Then, we have
\begin{equation}
\begin{split}
-C^\ast(t_{0}^{-1}) \int_M\varphi_\delta^2f\frac{\omega^n}{n!}&\leq \int_M\varphi_\delta^2f\Delta f\frac{\omega^n}{n!}\\
&=\int_M div(\varphi_\delta^2f\nabla f)\frac{\omega^n}{n!}-\int_M\nabla(\varphi_\delta^2f)\cdot\nabla f\frac{\omega^n}{n!}\\
&=-\int_M|\nabla(\varphi_\delta f)|_{\omega}^2\frac{\omega^n}{n!}+\int_M|\nabla\varphi_\delta|_{\omega}^2 f^2\frac{\omega^n}{n!}
\end{split}
\end{equation}
for $t\geq t_{0}>0$. By (\ref{cut}) and (\ref{Lp02}), we obtain
\begin{equation}
\begin{split}
 & \int_{M\setminus\Sigma}|\nabla f|_{\omega}^2\frac{\omega^n}{n!}=\lim_{\delta \rightarrow 0}\int_{M\setminus B_{2\delta}(\Sigma)}|\nabla f|_{\omega}^2\frac{\omega^n}{n!}\\
\leq &\lim_{\delta \rightarrow 0}\int_M|\nabla(\varphi_\delta f)|_{\omega}^2\frac{\omega^n}{n!}\\
\leq & \lim_{\delta \rightarrow 0} \int_M C^\ast(t_{0}^{-1})\varphi_\delta^2f+|\nabla\varphi_\delta|_{\omega}^2f^2\frac{\omega^n}{n!}\\
\leq & C^\ast(t_{0}^{-1})\cdot \hat{C}_{b}\\
\end{split}
\end{equation}
for $t\geq t_{0}>0$. This implies
 $f\in W^{1,2}(M, \omega)$ and $f$ satisfies the elliptic inequality  $\Delta f\geq -C^\ast(t_{0}^{-1})$ globally on $M$ in weakly sense for $t\geq t_{0}>0$. By the standard elliptic estimate (see Theorem 8.17 in \cite{GT}), we can show that $f\in L^{\infty}(M)$ for all $t\geq t_0>0$, and the $L^{\infty}$-norm depending on $C^\ast(t_{0}^{-1})$, the $L^{b}$-norm (i.e. $\hat{C}_{b}$) and the geometry of $(M, \omega)$, i.e. we have the following proposition.

 \medskip

\begin{proposition}\label{prop 3.4}
Along the heat flow (\ref{SSS1}), there exists a positive constant $\hat{C}_{\phi}$  depending only on $r, n, t_{0}^{-1}, C_{\phi}$ and the geometry of $(M, \omega)$ such that
\begin{equation}
\sup_{M\setminus \Sigma }|\phi|^2_{H(t), \omega }\leq \hat{C}_{\phi}
\end{equation}
for all $t\geq t_{0}>0$.
\end{proposition}

 \medskip

Recalling the Chern-Weil formula in \cite{Si} (Proposition 3.4) and using Fatou's lemma, we have
\begin{equation}\label{CW22}
\begin{split}
&4\pi^{2}\int_{M} (2c_{2}(\mathcal{E})-c_{1}(\mathcal{E})\wedge c_{1}(\mathcal{E}))\wedge\frac{\omega^{n-2}}{(n-2)!}\\
=&\lim_{\epsilon \rightarrow 0}4\pi^{2}\int_{\tilde{M}} (2c_{2}(E)-c_{1}(E)\wedge c_{1}(E))\wedge\frac{\omega_{\epsilon}^{n-2}}{(n-2)!}\\
=&\lim_{\epsilon \rightarrow 0}\int_{\tilde{M}}\tr (F_{H_{\epsilon}(t), \phi }\wedge F_{H_{\epsilon}(t), \phi })\wedge \frac{\omega_{\epsilon }^{n-2}}{(n-2)!}\\
=&\lim_{\epsilon \rightarrow 0}\int_{\tilde{M}}(|F_{H_{\epsilon}(t), \phi }|_{H_{\epsilon }(t), \omega_\epsilon}^{2}-|\Lambda_{\omega_{\epsilon}} F_{H_{\epsilon}(t), \phi }|_{H_{\epsilon }(t)}^{2}) \frac{\omega_{\epsilon}^{n}}{n!}\\
\geq&  \int_{M\setminus \Sigma}(|F_{H(t), \phi }|_{H(t), \omega}^{2}-|\sqrt{-1}\Lambda_{\omega} F_{H(t), \phi}|_{H(t)}^{2} )\frac{\omega^{n}}{n!}\\
\end{split}
\end{equation}
for $t>0$. Here,  over a non-projective compact complex manifold, the  Chern classes of a coherent sheaf can be defined by the classes of Atiyah-Hirzenbruch (\cite{AH}, see \cite{Gr} for details). The $L^{\infty}$ estimate of $|\phi|^2_{H(t), \omega }$, (\ref{H005}) and the above inequality imply that $|F_{H(t)}|_{H(t), \omega}$ is square integrable and $|\Lambda_{\omega } F_{H(t)}|_{H(t)}$ is uniformly bounded, i.e. we have the following corollary.

\medskip

\begin{corollary}\label{coro 3.5}
Let $H(t)$ be a solution of  the heat flow (\ref{SSS1}), then $H(t)$ must be an admissible Hermitian metric on $\mathcal{E}$ for every $t>0$.
\end{corollary}

\medskip

\section{Approximate Hermitian-Einstein structure }
\setcounter{equation}{0}

Let $H_{\epsilon}(t)$ be the long time solution of (\ref{DDD1}) and $H(t)$ be the long time solution of (\ref{SSS1}). We set:
\begin{equation}
\exp{S(t)}=h(t)=\hat{H}^{-1}H(t),
\end{equation}
\begin{equation}
\exp{S(t_{1}, t_{2})}=h(t_{1}, t_{2})=H^{-1}(t_{1})H(t_{2}),
\end{equation}
\begin{equation}
\exp{S_{\epsilon}(t_{1}, t_{2})}=h_{\epsilon}(t_{1}, t_{2})=H_{\epsilon}^{-1}(t_{1})H_{\epsilon}(t_{2}).
\end{equation}
By Lemma 3.1 in \cite{Si}, we have
\begin{equation}\label{la02}
\Delta_{\omega_{\epsilon }}\log(\tr h+\tr h^{-1})\geq -2|\Lambda_{\omega_{\epsilon} }(F_{H, \phi })|_{H}-2|\Lambda_{\omega }(F_{K, \phi })|_{K},
\end{equation}
where $\exp{S}=h=K^{-1}H$. By the uniform lower bound of Green functions $G_{\epsilon}$ (\ref{H005}) and the inequalities (\ref{c02}) , we have
\begin{equation}\label{for4.5}
\|S_{\epsilon }(t_{1}, t_{2})\|_{L^{\infty}(\tilde{M}) }\leq C_{1}\|S_{\epsilon }(t_{1}, t_{2})\|_{L^{1}(\tilde{M}, \omega_{\epsilon })}+C_{2}(t_{0}^{-1})
\end{equation}
for $0<t_{0}\leq t_{1} \leq t_{2}$, where $C_{1}$ is a constant depending only on the rank $r$ and  $C_{2}(t_{0}^{-1})$ is a constant depending only on $C_{K}$, $C_{G}$ and $t_{0}^{-1}$.
By limiting, we also have
\begin{equation}\label{mean1}
\|S(t_{1}, t_{2})\|_{L^{\infty}( M\setminus \Sigma )}\leq C_{1}\|S(t_{1}, t_{2})\|_{L^{1}( M\setminus \Sigma, \omega )}+C_{2}(t_{0}^{-1})
\end{equation}
for $0<t_{0}\leq t_{1} \leq t_{2}$.
On the other hand, (\ref{c01}) and (\ref{c02}) imply that
\begin{equation}\label{L101}
\begin{split}
& r^{-\frac{1}{2}}\|S(t_{1}, t_{2})\|_{L^1 ( M\setminus \Sigma, \omega )}-\Vol(M, \omega )\log(2r)\\ \leq& \int_{t_{1}}^{t_{2}}\int_{M\setminus \Sigma }|\sqrt{-1} \Lambda_{\omega} F_{H(s), \phi }-\lambda\Id_{\mathcal{E}}|_{H(s)}\frac{\omega^{n}}{n!} ds\\
\leq&\ \hat{C}_{1}(t_2-t_1). \\
\end{split}
\end{equation}
So, we know that the metrics $H(t_{1})$ and $H(t_{2})$ are  mutually bounded each other on $\mathcal{E}|_{M\setminus \Sigma}$. $(\mathcal{E}|_{M\setminus \Sigma}, \phi )$ can be seen as a Higgs bundle on the non-compact K\"ahler manifold $(M\setminus \Sigma , \omega)$. Let's recall  Donaldson's functional  defined on the space $\mathscr{P}_{0}$ of Hermitian metrics on the Higgs bundle $(\mathcal{E}|_{M\setminus \Sigma}, \phi )$ (see Section 5 in \cite{Si} for details),
\begin{eqnarray}\label{7}
\mu_{\omega} (K, H) = \int_{M\setminus \Sigma } \tr (S \sqrt{-1}\Lambda_{\omega }F_{K, \phi })+ \langle \Psi (S) (D''_{\phi} S) , D''_{\phi} S \rangle_{K}\frac{\omega^{n}}{n!},
\end{eqnarray}
where $\Psi (x, y)= (x-y)^{-2}(e^{y-x }-(y-x)-1)$, $\exp{S}=K^{-1}H$.
Since we have known  that $|\Lambda_{\omega } F_{H(t), \phi }|_{H(t)}$ is uniformly bounded for $t\geq t_{0}>0$, it is easy to see that $H(t)$ (for every $t>0$)  belongs to the definition space $\mathscr{P}_{0}$. By Lemma 7.1 in \cite{Si}, we have
a formula for the derivative with respect
to $t$ of Donaldson's functional,
\begin{eqnarray}\label{F5}
\frac{d}{dt}\mu (H(t_{1}), H(t)) = -2\int_{M\setminus \Sigma }|\Phi (H(t), \phi )|_{ H(t)}^{2}\frac{\omega^{n}}{n!}.
\end{eqnarray}

\medskip

 \begin{proposition}\label{prop 4.1}
 Let $H(t)$ be the long time solution of (\ref{SSS1}). If the reflexive Higgs sheaf $(\mathcal{E}, \phi )$ is $\omega$-semi-stable, then
\begin{equation}\label{semi03}
\int_{M\setminus \Sigma}|\sqrt{-1}\Lambda_{\omega} F_{H(t), \phi}-\lambda \Id_{\mathcal{E}} |_{H(t)}^{2} \frac{\omega^{n}}{n!}\rightarrow 0,
\end{equation}
as $t\rightarrow +\infty$.
\end{proposition}

\medskip

{\bf Proof. } We prove (\ref{semi03}) by contradiction. If not, by the monotonicity of  $\|\Lambda_\omega (F_{H(t), \phi})-\lambda\Id\|_{L^{2}}$, we can suppose that
\begin{equation}\lim_{t\rightarrow +\infty }\int_M|\sqrt{-1} \Lambda_\omega F_{H(t), \phi}-\lambda\Id_{\mathcal{E}}|_{H(t)}^{2}\frac{\omega^n}{n!} = C^\ast> 0.\end{equation}
By (\ref{F5}), we have
\begin{equation}\label{M05}
\mu_\omega(H(t_{0}), H(t))=-\int_{t_{0}}^t\int_{M\setminus \Sigma }|\Lambda_\omega F_{H(s), \phi }-\lambda\Id_{\mathcal{E}}|_{H(s)}^2\frac{\omega^n}{n!}ds\leq -C^\ast(t-t_{0})
\end{equation}
for all $0<t_{0} \leq t$. Then it is clear that (\ref{L101}) implies
\begin{equation}\label{semi01}
\liminf_{t\rightarrow +\infty}\frac{-\mu_{\omega } (H(t_{0}), H(t))}{\|S(t_{0}, t)\|_{L^{1}(M\setminus \Sigma , \omega )}}\geq r^{-\frac{1}{2}}\frac{C^\ast}{\hat{C}_{1}}.
\end{equation}
 By the definition of Donaldson's functional (\ref{7}), we must have a sequence  $t_{i}\rightarrow +\infty$ such that
\begin{equation}\label{CM02}
\|S(1, t_{i})\|_{L^{1}(M\setminus \Sigma , \omega)}\rightarrow +\infty.
\end{equation}
On the other hand, it is easy to check that
\begin{equation}\label{s01}
|S(t_{1}, t_{3})|_{H(t_{1})}\leq r(|S(t_{1}, t_{2})|_{H(t_{1})}+|S(t_{2}, t_{3})|_{H(t_{2})})
\end{equation}
for all $0\leq t_{1} ,  t_{2} , t_{3}$. Then, by (\ref{mean1}), we have
\begin{equation}\label{CM03}
\lim_{i\rightarrow \infty }\|S(t_{0}, t_{i})\|_{L^{1}(M\setminus \Sigma , \omega)}\rightarrow +\infty ,
\end{equation}
 and
\begin{equation}\label{CM01}
\begin{split}
&\|S(t_{0}, t)\|_{L^{\infty}(M\setminus \Sigma )}\leq  r\|S(1, t)\|_{L^{\infty}(M\setminus \Sigma )}+r \|S(t_{0}, 1)\|_{L^{\infty}(M\setminus \Sigma )}\\
 \leq &\ r^{2}C_{3}(\|S(t_{0}, t)\|_{L^{1}}+\|S(t_{0}, 1)\|_{L^{1}})+r\|S(t_{0}, 1)\|_{L^{\infty}(M\setminus \Sigma )}+rC_{4}\\
 \end{split}
\end{equation}
 for all $0<t_{0}\leq t$, where $C_{3}$ and $C_{4}$ are uniform constants depending only on $r$, $C_{K}$ and $C_{G}$.

Set $u_{i}(t_{0})=\|S(t_{0}, t_{i})\|_{L^{1}}^{-1}S(t_{0}, t_{i})\in S_{H(t_{0})}(\mathcal{E}|_{M\setminus \Sigma })$, where $S_{H(t_{0})}(\mathcal{E}|_{M\setminus \Sigma })=\{\eta \in \Omega^{0}(M\setminus \Sigma , \End(\mathcal{E}|_{M\setminus \Sigma }))| \quad \eta ^{\ast H(t_{0})}=\eta \}$, then $\|u_{i}(t_{0})\|_{L^{1}}=1$.
By (\ref{1}) and (\ref{for4.5}), we have
\begin{equation}
\int_{M\setminus\Sigma}\tr S(t_{0}, t_i) \frac{\omega^{n}}{n!}=0,
\end{equation}
so
\begin{equation}\int_{M\setminus\Sigma}\tr u_{i}(t_{0}) \frac{\omega^{n}}{n!}=0.\end{equation}
   By the inequalities (\ref{semi01}), (\ref{CM02}), (\ref{CM01}),   and the Lemma 5.4 in \cite{Si}, we can see that,
by choosing a subsequence which we also denote by $u_i(t_0)$, we have $u_{i}(t_{0})\rightarrow u_{\infty}(t_{0})$   weakly  in $L_{1}^{2}$, where the limit $u_{\infty}(t_{0})$ satisfies: $\|u_{\infty}(t_{0})\|_{L^{1}}=1$, $\int_{M}\tr(u_{\infty}(t_{0}))\frac{\omega^{n} }{n!}=0$ and
\begin{equation}\label{t01}
 \|u_{\infty}(t_{0})\|_{L^{\infty}}\leq r^{2}C_{3}.
 \end{equation}
Furthermore, if $\Upsilon : R\times R \rightarrow R$ is a positive smooth function such that $\Upsilon (\lambda_{1}, \lambda_{2})< (\lambda_{1}- \lambda_{2})^{-1}$ whenever $\lambda_{1}>\lambda_{2}$, then
 \begin{equation}\label{t02}
 \begin{split}
 &\int_{M\setminus \Sigma}\tr (u_{\infty}(t_{0})\sqrt{-1}\Lambda_{\omega }(F_{H(t_{0}), \phi })) + \langle \Upsilon  (u_{\infty}(t_{0}))(\overline{\partial }_{\phi }u_{\infty}(t_{0})), \overline{\partial }_{\phi }u_{\infty}(t_{0}) \rangle_{H(t_{0})}\frac{\omega^{n} }{n!}\\&\leq -r^{-\frac{1}{2}}\frac{C^\ast}{\hat{C}_{1}}.\\
 \end{split}
 \end{equation}

 \medskip

Since $ \|u_{\infty}(t_{0})\|_{L^{\infty}}$ and $\|\Lambda_{\omega }(F_{H(t_{0}), \phi })\|_{L^{1}}$ are uniformly bounded (independent of $t_{0}$), (\ref{t02}) implies that: there exists a uniform constant $\check{C}$ independent of $t_{0}$ such that
\begin{equation}\label{t03}
 \int_{M\setminus \Sigma } |\overline{\partial }_{\phi}u_{\infty}(t_{0})|_{H(t_{0})}^{2}\frac{\omega^{n} }{n!}\leq \check{C}.
 \end{equation}

From Lemma \ref{lem 2.2}, we see that $\hat{H}$ and $H(t_{0})$ are locally mutually bounded each other. By choosing a subsequence, we have $u_{\infty}(t_{0}) \rightarrow u_{\infty}$ weakly in local $L_{1}^{2}$ outside $\Sigma $ as $t_{0}\rightarrow 0$, where $u_{\infty}$ satisfies
\begin{equation}\int_{M}\tr(u_{\infty})\frac{\omega^{n} }{n!}=0, \quad and \quad \|u_{\infty}\|_{L^{1}}=1.\end{equation}

Since $|\sqrt{-1} \Lambda_{\omega_\epsilon}F_{H_\epsilon(t), \phi }|_{H_\epsilon(t)}\in L^{\infty}$ for $t>0$, by the uniform upper bound of the heat kernels (\ref{kernel01}), we have
\begin{equation}\label{FL1}
\begin{split}
&\int_{B_{\omega_1}(\delta)\setminus \Sigma}|\sqrt{-1} \Lambda_{\omega}F_{H(t), \phi }|_{H(t)}\frac{\omega^n}{n!}\\
=&\lim_{\epsilon\rightarrow 0}\int_{B_{\omega_1}(\delta)}|\sqrt{-1} \Lambda_{\omega_\epsilon}F_{H_\epsilon(t), \phi }|_{H_\epsilon(t)}\frac{\omega_\epsilon^n}{n!}\\
\leq &\lim_{\epsilon\rightarrow 0}\int_{B_{\omega_1}(\delta)}\int_{\tilde{M}}K_{\epsilon}(x, y, t)|\sqrt{-1} \Lambda_{\omega_\epsilon}F_{\hat{H}, \phi }|_{\hat{H}}(y)\frac{\omega_\epsilon^n(y)}{n!}\cdot\frac{\omega_\epsilon^n(x)}{n!}\\
=&\lim_{\epsilon\rightarrow 0}\int_{B_{\omega_1}(\delta)}\Big(\big(\int_{B_{\omega_1}(2\delta)}+\int_{\tilde{M}\setminus B_{\omega_1}(2\delta)}\big)K_{\epsilon}(x, y, t)|\sqrt{-1} \Lambda_{\omega_\epsilon}F_{\hat{H}, \phi}|_{\hat{H}}(y)\frac{\omega_\epsilon^n(y)}{n!}\Big)\frac{\omega_\epsilon^n(x)}{n!}\\
\leq &\lim_{\epsilon\rightarrow 0}\int_{\tilde{M}}\int_{B_{\omega_1}(2\delta)}K_{\epsilon}(x, y, t)|\sqrt{-1} \Lambda_{\omega_\epsilon}F_{\hat{H}, \phi}|_{\hat{H}}(y)\frac{\omega_\epsilon^n(y)}{n!}\cdot\frac{\omega_\epsilon^n(x)}{n!}\\
&+\int_{B_{\omega_1}(\delta)}\Big(\int_{\tilde{M}\setminus B_{\omega_1}(2\delta)}C_K(\tau)t^{-n}\exp\big(-\frac{d_{\omega_\epsilon}(x, y)}{(4+\tau)t}\big)|\sqrt{-1} \Lambda_{\omega_\epsilon}F_{\hat{H}, \phi}|_{\hat{H}}(y)\frac{\omega_\epsilon^n(y)}{n!}\Big)\frac{\omega_\epsilon^n(x)}{n!}\\
\leq &\int_{B_{\omega_1}(2\delta)\setminus \Sigma}|\sqrt{-1} \Lambda_{\omega}F_{\hat{H}, \phi}|_{\hat{H}}\frac{\omega^n}{n!}\\
&+C_K(\tau)t^{-n}\exp\big(-\frac{a(\delta)}{(4+\tau)t}\big)\Vol_{\omega_1}(B_{\omega_1}(\delta))\int_{M}|\sqrt{-1} \Lambda_{\omega}F_{\hat{H}, \phi}|_{\hat{H}}\frac{\omega^n}{n!}.\\
\end{split}
\end{equation}
By (\ref{FL1}) and the uniform bound of $\|u_{\infty}(t_{0})\|_{L^{\infty}}$, we have
\begin{equation}\label{t04}
 \lim_{t_{0}\rightarrow 0 }\int_{M}\tr (u_{\infty}(t_{0})\sqrt{-1}\Lambda_{\omega }F_{H(t_{0}), \phi }) \frac{\omega^{n} }{n!}= \int_{M}\tr (u_{\infty}\sqrt{-1}\Lambda_{\omega }F_{\hat{H}, \phi}) \frac{\omega^{n} }{n!}.
 \end{equation}

Let's denote
\begin{eqnarray}\label{7.2}
S_{\hat{H}}(\mathcal{E}|_{M\setminus \Sigma })=\{\eta \in \Omega^{0}(M\setminus \Sigma , \End(\mathcal{E}|_{M\setminus \Sigma }))| \quad \eta ^{\ast \hat{H}}=\eta \}.
\end{eqnarray}
and
\begin{equation}\hat{u}_{\infty}(t_{0})=(h(t_{0}))^{\frac{1}{2}}\cdot u_{\infty}(t_{0}) \cdot (h(t_{0}))^{-\frac{1}{2}}.\end{equation} It is easy to check that: $\hat{u}_{\infty}(t_{0})\in S_{\hat{H}}(\mathcal{E}|_{M\setminus \Sigma })$ and $|\hat{u}_{\infty}(t_{0})|_{\hat{H}}=|u_{\infty}(t_{0})|_{H(t_{0})}$. Furthermore, we have:

\medskip

\begin{lemma}\label{lem 4.2}
For any compact domain $\Omega \subset M\setminus \Sigma$ and any positive smooth function $\Upsilon : R\times R \rightarrow R$, we have
\begin{equation}\label{t05}
\lim_{t_{0}\rightarrow 0}\int_{\Omega} |\langle\Upsilon  (u_{\infty}(t_{0}))(\overline{\partial }_{\phi}u_{\infty}(t_{0})), \overline{\partial }_{\phi}u_{\infty}(t_{0})\rangle_{H(t_{0})}-\langle\Upsilon  (\hat{u}_{\infty}(t_{0}))(\overline{\partial }_{\phi}\hat{u}_{\infty}(t_{0})), \overline{\partial }_{\phi}\hat{u}_{\infty}(t_{0})\rangle_{\hat{H}}|\frac{\omega^{n} }{n!}=0.
\end{equation}
\end{lemma}

\medskip

{\bf Proof. } At each point $x$ on $\Omega $, we choose
a unitary basis $\{e_{i}\}_{i=1}^{r}$  with respect to the metric
$H(t_{0})$, such that $u_{\infty}(t_{0}) (e_{i})= \lambda_{i} e_{i}$. Then, $\{\hat{e}_{i}=(h(t_{0}))^{\frac{1}{2}}e_{i}\}$ is a unitary basis with respect to the metric $\hat{H}$ and
$\hat{u}_{\infty}(t_{0}) (\hat{e}_{i})= \lambda_{i} \hat{e}_{i}$. Set:
 \begin{equation}
\overline{\partial }_{\phi}u_{\infty}(t_{0})(e_{i})=(\overline{\partial }_{\phi}u_{\infty}(t_{0}))_{i}^{j}e_{j}, \quad \overline{\partial }_{\phi}\hat{u}_{\infty}(t_{0})(\hat{e}_{i})=(\overline{\partial }_{\phi}\hat{u}_{\infty}(t_{0}))_{i}^{j}\hat{e}_{j},
\end{equation}
then
\begin{equation}
|\overline{\partial }_{\phi}u_{\infty}(t_{0})|_{H(t_{0}), \omega}^{2}=\sum_{i, j=1}^{r}\langle(\overline{\partial }_{\phi}u_{\infty}(t_{0}))_{i}^{j}, (\overline{\partial }_{\phi}u_{\infty}(t_{0}))_{i}^{j}\rangle_{\omega},
\end{equation}
\begin{equation}\label{z1}
 \langle\Upsilon  (u_{\infty}(t_{0}))(\overline{\partial }_{\phi}u_{\infty}(t_{0})), \overline{\partial }_{\phi}u_{\infty}(t_{0})\rangle_{H(t_{0})}
=\sum_{i, j=1}^{r}\langle\Upsilon (\lambda_{i}, \lambda_{j})(\overline{\partial }_{\phi}u_{\infty}(t_{0}))_{i}^{j}, (\overline{\partial }_{\phi}u_{\infty}(t_{0}))_{i}^{j}\rangle_{\omega},
\end{equation}
\begin{equation}
\Upsilon  (\hat{u}_{\infty}(t_{0}))(\overline{\partial }_{\phi}\hat{u}_{\infty}(t_{0}))(\hat{e}_{i})=\sum_{j=1}^{r}\Upsilon (\lambda_{i}, \lambda_{j})(\overline{\partial }_{\phi}\hat{u}_{\infty}(t_{0}))_{i}^{j}\hat{e}_{j},
\end{equation}
and
\begin{equation}\label{z2}
 \langle\Upsilon  (\hat{u}_{\infty}(t_{0}))(\overline{\partial }_{\phi}\hat{u}_{\infty}(t_{0})), \overline{\partial }_{\phi}\hat{u}_{\infty}(t_{0})\rangle_{\hat{H}}
=\sum_{i, j=1}^{r}\langle\Upsilon (\lambda_{i}, \lambda_{j})(\overline{\partial }_{\phi}\hat{u}_{\infty}(t_{0}))_{i}^{j}, (\overline{\partial }_{\phi}\hat{u}_{\infty}(t_{0}))_{i}^{j}\rangle_{\omega}.
\end{equation}
By the definition, we have
\begin{equation}
\begin{split}
\overline{\partial }_{\phi}\hat{u}_{\infty}(t_{0})=&\ (h(t_{0}))^{\frac{1}{2}}\circ \overline{\partial }_{\phi}u_{\infty}(t_{0}) \circ (h(t_{0}))^{-\frac{1}{2}} +\overline{\partial }_{\phi}(h(t_{0}))^{\frac{1}{2}} \circ
u_{\infty}(t_{0}) \circ (h(t_{0}))^{-\frac{1}{2}}\\
& -(h(t_{0}))^{\frac{1}{2}} \circ u_{\infty}(t_{0}) \circ (h(t_{0}))^{-\frac{1}{2}}\circ \overline{\partial }_{\phi}(h(t_{0}))^{\frac{1}{2}}\circ (h(t_{0}))^{-\frac{1}{2}}\\
=&\ (h(t_{0}))^{\frac{1}{2}}\circ \overline{\partial }_{\phi}u_{\infty}(t_{0}) \circ (h(t_{0}))^{-\frac{1}{2}} +\overline{\partial }_{\phi}(h(t_{0}))^{\frac{1}{2}} \circ (h(t_{0}))^{-\frac{1}{2}}
\hat{u}_{\infty}(t_{0})\\
& -\hat{u}_{\infty}(t_{0})\circ \overline{\partial }_{\phi}(h(t_{0}))^{\frac{1}{2}}\circ (h(t_{0}))^{-\frac{1}{2}},\\
\end{split}
\end{equation}
and
\begin{equation}\label{z3}
(\overline{\partial }_{\phi}\hat{u}_{\infty}(t_{0}))_{i}^{j} =(\overline{\partial }_{\phi}u_{\infty}(t_{0}))_{i}^{j} +(\lambda_{i}-\lambda_{j})\{\overline{\partial }_{\phi}(h(t_{0})^{\frac{1}{2}} \circ (h(t_{0}))^{-\frac{1}{2}}\}_{i}^{j},
\end{equation}
where $\overline{\partial }_{\phi}(h(t_{0})^{\frac{1}{2}}\circ (h(t_{0})^{-\frac{1}{2}})(\hat{e}_{i})=(\overline{\partial }_{\phi}(h(t_{0})^{\frac{1}{2}} \circ (h(t_{0})^{-\frac{1}{2}})_{i}^{j}\hat{e}_{j}$.
By (\ref{t01}), (\ref{z1}), (\ref{z2}) and (\ref{z3}), we have
\begin{equation}\label{z4}
\begin{split}
& |\langle\Upsilon  (\hat{u}_{\infty}(t_{0}))(\overline{\partial }_{\phi}\hat{u}_{\infty}(t_{0})), \overline{\partial }_{\phi}\hat{u}_{\infty}(t_{0})\rangle_{\hat{H}}-\langle\Upsilon  (u_{\infty}(t_{0}))(\overline{\partial }_{\phi}u_{\infty}(t_{0})), \overline{\partial }_{\phi}u_{\infty}(t_{0})\rangle_{H(t_{0})}|\\
&\leq 8(r^{2}C_{3})^{2}(B^{\ast}(\Upsilon) ) (|\overline{\partial }_{\phi}u_{\infty}(t_{0})|_{H(t_{0})}|\overline{\partial }_{\phi}(h(t_{0})^{\frac{1}{2}} \circ (h(t_{0}))^{-\frac{1}{2}}|_{\hat{H}}+|\overline{\partial }_{\phi}(h(t_{0})^{\frac{1}{2}} \circ (h(t_{0}))^{-\frac{1}{2}}|_{\hat{H}}^{2}),\\
\end{split}
\end{equation}
where $B^{\ast}(\Upsilon)=\max_{[-r^{2}C_{3}, r^{2}C_{3}]^{2}}\Upsilon$.
Since $H(t)$ are smooth on $M\setminus \Sigma \times [0, 1]$ and $h(t)\rightarrow \Id_{\mathcal{E}}$ locally in $C^{\infty}$-topology as $t\rightarrow 0$, it is easy to check that
\begin{equation}\label{z5}
\sup_{x\in \Omega}(|(h(t_{0}))^{-\frac{1}{2}}\overline{\partial}_{\phi }(h(t_{0}))^{\frac{1}{2}}|_{\hat{H}, \omega}+|\overline{\partial}_{\phi }(h(t_{0}))^{\frac{1}{2}}(h(t_{0}))^{-\frac{1}{2}}|_{\hat{H}, \omega})\leq C_{\Omega}(t_{0}),
\end{equation}
where $C_{\Omega}(t_{0})\rightarrow 0$ as $t_{0}\rightarrow 0$. On the other hand, $|\overline{\partial }_{\phi }u_{\infty}(t_{0})|_{H(t_{0}), \omega}$ are uniform bounded in $L^{2}$, so (\ref{z4}) and (\ref{z5}) imply (\ref{t05}).

\hfill $\Box$ \\

By (\ref{t02}), (\ref{t04}) and (\ref{t05}), we have that given any compact domain $\Omega \subset M\setminus \Sigma$ and any positive number $\tilde{\epsilon}>0$,
\begin{equation}\label{t07}
 \int_{M\setminus\Sigma}\tr (u_{\infty}\sqrt{-1}\Lambda_{\omega }F_{\hat{H}}, \phi )\frac{\omega^{n} }{n!}+\int_{\Omega} \langle\Upsilon  (\hat{u}_{\infty}(t_{0}))(\overline{\partial }_{\phi }\hat{u}_{\infty}(t_{0})), \overline{\partial }_{\phi }\hat{u}_{\infty}(t_{0})\rangle_{\hat{H}}\frac{\omega^{n} }{n!}\leq -r^{-\frac{1}{2}}\frac{C^\ast}{\hat{C}_{1}}+\tilde{\epsilon}
 \end{equation}
for small $t_{0}$. As we know that $\hat{u}_{\infty}(t_{0}) \rightarrow u_{\infty}$ in $L^{2}(\Omega )$, $|\hat{u}_{\infty}(t_{0})|_{\hat{H}}$ is uniformly bounded in $L^{\infty}$ and $|\overline{\partial }_{\phi}\hat{u}_{\infty}(t_{0})|_{\hat{H}, \omega}$ is uniformly bounded in $L^{2}(\Omega )$. By the same argument as that in Simpson's paper (Lemma 5.4 in \cite{Si}), we have
\begin{equation}\label{t07}
 \int_{M\setminus\Sigma}\tr (u_{\infty}\sqrt{-1}\Lambda_{\omega }F_{\hat{H}, \phi })\frac{\omega^{n} }{n!}+\|\Upsilon^{\frac{1}{2}}  (u_{\infty})(\overline{\partial }_{\phi}u_{\infty})\|_{L^{q}(\Omega)}^{2}\leq -r^{-\frac{1}{2}}\frac{C^\ast}{\hat{C}_{1}}+2\tilde{\epsilon}
 \end{equation}
for any $q<2$ and any $\tilde{\epsilon}$. Since $\tilde{\epsilon}$, $q<2$ and $\Omega$ are arbitrary,  we get
\begin{equation}\label{semi02}
 \int_{M\setminus\Sigma}\tr (u_{\infty}\sqrt{-1}\Lambda_{\omega }F_{\hat{H}, \phi}) + \langle\Upsilon  (u_{\infty})(\overline{\partial }_{\phi}u_{\infty}), \overline{\partial }_{\phi}u_{\infty}\rangle_{\hat{H}}\frac{\omega^{n} }{n!}\leq -r^{-\frac{1}{2}}\frac{C^\ast}{\hat{C}_{1}}.
 \end{equation}

\medskip

By the above inequality and the Lemma 5.5 in \cite{Si}, we can see that
the eigenvalues of $u_{\infty}$ are  constant almost everywhere.
Let $\lambda_{1} < \dots <\lambda _{l}$ denote the distinct eigenvalue of $u_{\infty}$.  Since $\int_{M} \tr u_{\infty} \frac{\omega^{n} }{n!}=0$ and $\|u_{\infty}\|_{L^{1}}=1$, we must have $l\geq 2$.
For any $1\leq \alpha <l $, define function $P_{\alpha } : R\rightarrow R$ such that
\begin{eqnarray}
P_{\alpha }=\left \{\begin{array}{cll} 1, & x\leq \lambda_{\alpha },\\
0,
 & x\geq \lambda_{\alpha +1}.\\
\end{array}\right.
\end{eqnarray}
Set $\pi_{\alpha }=P_{\alpha } (u_{\infty})$, Simpson (p887 in \cite{Si}) proved that:

\medskip

(1) $\pi_{\alpha } \in L_{1}^{2}(M\setminus \Sigma , \omega , \hat{H})$;

(2)  $\pi_{\alpha }^{2}=\pi_{\alpha }=\pi_{\alpha }^{\ast \hat{H}}$;

(3) $(\Id_{\mathcal{E}} -\pi_{\alpha }) \bar{\partial }\pi_{\alpha } =0$;

(4) $(\Id_{\mathcal{E}} -\pi_{\alpha }) [\phi , \pi_{\alpha }]=0$.

\medskip

By Uhlenbeck and Yau's regularity statement of $L_{1}^{2}$-subbundle (\cite{UY}), $\pi_{\alpha }$ represent a saturated coherent  Higgs sub-sheaf $E_{\alpha }$ of $(\mathcal{E}, \phi )$ on the open set $M\setminus \Sigma $. Since the singularity set $\Sigma $ is co-dimension at least $3$, by Siu's extension theorem (\cite{Siu1}), we know that  $E_{\alpha }$ admits a coherent analytic extension $\tilde{E}_{\alpha }$. By Serre's result (\cite{Se}), we get the direct image $i_{\ast }E_{\alpha }$ under the inclusion $i: M\setminus \Sigma \rightarrow M$ is coherent. So, every $E_{\alpha }$  can be extended to the whole $M$ as a saturated coherent Higgs sub-sheaf of $(\mathcal{E}, \phi )$, which will also be denoted by $E_{\alpha }$ for simplicity. By the Chern-Weil formula (\ref{CW1}) (Proposition 4.1 in \cite{Br}) and the above condition (4), we have
\begin{equation}
\begin{split}
\deg_{\omega} (E_{\alpha })&=\int_{M\setminus \Sigma } \tr (\pi_{\alpha } \sqrt{-1}\Lambda _{\omega }F_{\hat{H} }) -|\overline{\partial } \pi_{\alpha }|_{\hat{H}, \omega}^{2} \frac{\omega^{n}}{n!}\\
&=\int_{M\setminus \Sigma } \tr (\pi_{\alpha } \sqrt{-1}\Lambda _{\omega }F_{\hat{H}, \phi }) -|D''_{\phi } \pi_{\alpha }|_{K, \omega}^{2} \frac{\omega^{n}}{n!}.
\end{split}
\end{equation}
Set
\begin{eqnarray}
\nu =\lambda _{l} \deg_\omega (\mathcal{E}) -\sum_{\alpha =1} ^{l-1} (\lambda_{\alpha +1 } -\lambda_{\alpha }) \deg_\omega (E_{\alpha}).
\end{eqnarray}
Since $u_{\infty }=\lambda _{l} \Id_{\mathcal{E}}  -\sum_{\alpha =1} ^{l-1} (\lambda_{\alpha +1} -\lambda_{\alpha })\pi_{\alpha }$ and $\int_{M\setminus \Sigma }\tr u_{\infty }\frac{\omega^{n}}{n!}=0$, we have
\begin{eqnarray}
\lambda _{l} \rank (\mathcal{E}) -\sum_{\alpha =1} ^{l-1}(\lambda_{\alpha +1} -\lambda_{\alpha }) \rank (E_{\alpha })=0,
\end{eqnarray}
then
\begin{eqnarray}\label{3}
\nu =\sum_{\alpha =1} ^{l-1} (\lambda_{\alpha +1} -\lambda_{\alpha }) \rank (E_{\alpha }) (\frac{\deg_\omega (\mathcal{E})}{\rank (\mathcal{E})}-\frac{\deg_\omega (E_{\alpha })}{\rank (E_{\alpha })}).
\end{eqnarray}

By the argument similar to the one used in Simpson's paper (P888 in \cite{Si}) and the inequality (\ref{semi02}), we have
\begin{equation}
\begin{split}
\nu =& \int_{M}\tr (u_{\infty }\sqrt{-1}\Lambda _{\omega }F_{\hat{H}, \phi })\\
& +\langle \sum_{\alpha =1} ^{l-1} (\lambda_{\alpha +1}-\lambda_{\alpha })(dP_{\alpha })^{2}(u_{\infty }) (D''_{\phi} u_{\infty}) , D''_{\phi} u_{\infty}\rangle _{\hat{H}}\\
\leq&  -r^{-\frac{1}{2}}\frac{C^\ast}{\hat{C}_{1}}.\\
\end{split}
\end{equation}
On the other hand,  (\ref{3}) and the semi-stability imply $\nu \geq 0$, so we get a contradiction.

\hfill $\Box$ \\

{\bf Proof of Theorem \ref{thm 1.1}}\ By (\ref{H00001}), we have
\begin{equation}
\sup_{x\in M\setminus \Sigma}|\sqrt{-1} \Lambda_\omega (F_{H(t+1), \phi })-\lambda\Id_{\mathcal{E}}|_{H(t+1)}^2(x)\leq C_{K}\int_{M\setminus \Sigma}|\sqrt{-1} \Lambda_\omega (F_{H(t), \phi})-\lambda\Id_{\mathcal{E}}|_{H(t)}^2\frac{\omega^n}{n!}.
\end{equation}
If the reflexive Higgs sheaf $(\mathcal{E}, \phi )$ is $\omega $-semi-stable, (\ref{semi03}) implies
\begin{equation}
\sup_{x\in M\setminus \Sigma}|\sqrt{-1} \Lambda_\omega (F_{H(t), \phi })-\lambda\Id_{\mathcal{E}}|_{H(t+1)}^2(x)\rightarrow 0,
\end{equation}
as $t\rightarrow +\infty$. By corollary \ref{coro 3.5}, we know that every $H(t)$ is an admissible Hermitian metric. Then we get an approximate Hermitian-Einstein structure on a semi-stable reflexive Higgs sheaf.

By choosing a subsequence $\epsilon \rightarrow 0$, we have $H_{\epsilon}(t)$ converge to $H(t)$ in local $C^{\infty}$-topology. Applying Fatou's lemma we obtain
\begin{equation}
\begin{split}
&4\pi^{2}\int_{M} (2c_{2}(\mathcal{E})-\frac{r-1}{r}c_{1}(\mathcal{E})\wedge c_{1}(\mathcal{E}))\wedge\frac{\omega^{n-2}}{(n-2)!}\\
=&\lim_{\epsilon \rightarrow 0}4\pi^{2}\int_{\tilde{M}} (2c_{2}(E)-\frac{r-1}{r}c_{1}(E)\wedge c_{1}(E))\wedge\frac{\omega_{\epsilon}^{n-2}}{(n-2)!}\\
=&\lim_{\epsilon \rightarrow 0}\int_{\tilde{M}}\tr (F_{H_{\epsilon}(t), \phi }^{\bot}\wedge F_{H_{\epsilon}(t), \phi }^{\bot})\wedge \frac{\omega_{\epsilon }^{n-2}}{(n-2)!}\\
=&\lim_{\epsilon \rightarrow 0}\int_{\tilde{M}}|F_{H_{\epsilon}(t), \phi }^{\bot }|_{H_{\epsilon }(t), \omega_\epsilon}^{2}-|\Lambda_{\omega_{\epsilon}} F_{H_{\epsilon}(t), \phi }^{\bot}|_{H_{\epsilon }(t)}^{2} \frac{\omega_{\epsilon}^{n}}{n!}\\
\geq&  \int_{M\setminus \Sigma}|F_{H(t), \phi }^{\bot }|_{H(t), \omega}^{2}\frac{\omega^{n}}{n!}\\ & -\int_{M\setminus \Sigma}|\sqrt{-1}\Lambda_{\omega} F_{H(t), \phi}-\lambda \Id_{\mathcal{E}} -\frac{1}{r} \tr (\sqrt{-1}\Lambda_{\omega} F_{H(t), \phi }-\lambda \Id_{\mathcal{E}})\Id_{\mathcal{E}}|_{H(t)}^{2} \frac{\omega^{n}}{n!}\\
\end{split}
\end{equation}
for $t>0$, where $F_{H, \phi }^{\bot}$ is the trace free part of $F_{H, \phi }$. Let $t\rightarrow +\infty$, then (\ref{semi03}) implies the following Bogomolov type inequality
\begin{equation}\label{Bog}
\int_{M} (2c_{2}(\mathcal{E})-\frac{r-1}{r}c_{1}(\mathcal{E})\wedge c_{1}(\mathcal{E}))\wedge\frac{\omega^{n-2}}{(n-2)!}\geq 0.
\end{equation}

Now we prove that the existence of an approximate Hermitian-Einstein structure implies the semistability of $(\mathcal{E}, \phi )$.
Let $s$ be a $\theta $-invariant holomorphic section of a reflexive Higgs sheaf $(\mathcal{G}, \theta )$ on a compact K\"ahler manifold $(M, \omega )$, i.e. there exists a holomorphic $1$-form $\eta $ on $M\setminus \Sigma_{\mathcal{G}}$ such that $\theta (s)=\eta \otimes s$, where $\Sigma_{\mathcal{G}}$ is the singularity set of $\mathcal{G}$. Given a Hermitian metric $H$ on $\mathcal{G}$,
by computing, we have
\begin{equation}\label{W01}
\begin{split}
 &\sqrt{-1}\Lambda_{\omega } \langle s, -[\theta , \theta^{\ast H}]s\rangle_{H}\\
= & -\sqrt{-1}\Lambda_{\omega } \langle\theta^{\ast H}s,  \theta^{\ast H}s\rangle_{H}-\sqrt{-1}\Lambda_{\omega } \langle\theta s,  \theta s\rangle_{H}\\
= & -\sqrt{-1}\Lambda_{\omega } \langle\theta^{\ast H}s- \langle\theta^{\ast H}s,  s\rangle_{H}\frac{s}{|s|_{H}^{2}},  \theta^{\ast H}s - \langle\theta^{\ast H}s,  s\rangle_{H}\frac{s}{|s|_{H}^{2}}\rangle_{H}\\
  & -\sqrt{-1}\Lambda_{\omega } \langle \langle\theta^{\ast H}s,  s\rangle_{H}\frac{s}{|s|_{H}^{2}},   \langle\theta^{\ast H}s,  s\rangle_{H}\frac{s}{|s|_{H}^{2}}\rangle_{H}-\sqrt{-1}\Lambda_{\omega } \langle\phi s,  \phi s\rangle_{H}\\
= &\ |\theta^{\ast H}s- \langle\theta^{\ast H}s,  s\rangle_{H}\frac{s}{|s|_{H}^{2}}|_{H, \omega}^{2}\geq 0,\\
\end{split}
\end{equation}
where we have used $\theta (s)=\eta \otimes s$ in the third equality. Then, we have the following Weitzenb\"ock formula
\begin{equation}
\begin{split}
 \frac{1}{2}\Delta_{\omega } |s|_{H}^{2}=&\ \sqrt{-1}\Lambda_{\omega }\partial \overline{\partial}|s|_{H}^{2}\\
= &\ |D_{H}^{1, 0}s|_{H, \omega}^{2}+\sqrt{-1}\Lambda_{\omega }\langle s, F_{H}s\rangle_{H} \\
= &\ |D_{H}^{1, 0}s|_{H, \omega}^{2}-\langle s, \sqrt{-1}\Lambda_{\omega }F_{H, \theta }s\rangle_{H}-\sqrt{-1}\Lambda_{\omega } \langle s, [\theta , \theta^{\ast H}]s\rangle_{H} \\
\geq &\ |D_{H}^{1, 0}s|_{H, \omega}^{2}-\langle s, \sqrt{-1}\Lambda_{\omega }F_{H, \theta }s\rangle_{H}
\end{split}
\end{equation}
on $M\setminus \Sigma_{\mathcal{G}}$.

We suppose that the reflexive  Higgs sheaf $(\mathcal{G}, \theta )$ admits an approximate admissible Hermitian-Einstein structure, i.e.  for every positive $\delta
$, there is an admissible Hermitian metric $H_{\delta}$ such that
\begin{equation}
\sup _{x\in M\setminus \Sigma_{\mathcal{G}} } |\sqrt{-1}\Lambda_{\omega }F_{H_{\delta}, \theta }-\lambda(\mathcal{G}) \Id|_{H_{\delta}}(x)<\delta .
\end{equation}
If $\deg_{\omega}\mathcal{G}$ is negative, i.e. $\lambda(\mathcal{G})<0$, by choosing $\delta$ small enough, we have
\begin{equation}\label{la}
\Delta_{\omega } |s|_{H_{\delta}}^{2}
\geq 2 |D_{H}^{1, 0}s|_{H_{\delta}, \omega}^{2}-\lambda(\mathcal{G})|s|_{H_{\delta}}^{2}
\end{equation}
on  $M\setminus \Sigma_{\mathcal{G}}$.
Since every $H_{\delta }$ is admissible, by Theorem 2 in \cite{BS}, we know that $|s|_{H_{\delta}}\in L^{\infty}(M)$. Then,  the inequality (\ref{la}) can be extended globally to the compact manifold $M$. So, we must have
\begin{equation}
s\equiv 0.
\end{equation}

Assume that $(\mathcal{E}, \phi )$ admits an approximate Hermitian-Einstein structure and $\mathcal{F}$ is a saturated Higgs subsheaf of $(\mathcal{E}, \phi )$ with rank $p$. Let $\mathcal{G}=\wedge^{p}\mathcal{E}\otimes \det(\mathcal{F})^{-1}$, and $\theta $ be a Higgs filed naturally induced on $\mathcal{G}$ by the Higgs field $\phi $. One can check that $(\mathcal{G}, \theta )$ is also a reflexive Higgs sheaf which admits an approximate Hermitian-Einstein structure with constant
\begin{equation}
\lambda(\mathcal{G})=\frac{2p\pi}{\Vol(M, \omega )}(\mu_{\omega}(\mathcal{E})-\mu_{\omega}(\mathcal{F})).
\end{equation}
The inclusion $\mathcal{F}\hookrightarrow \mathcal{E}$ induces a morphism $\det(\mathcal{F})\rightarrow \wedge^{p}\mathcal{E} $ which can be seen as a nontrivial $\theta$-invariant holomorphic section of  $\mathcal{G} $. From above, we have $\lambda(\mathcal{G})\geq 0$, so the reflexive sheaf $(\mathcal{E}, \phi )$ is $\omega$-semistable.
 This completes the proof of Theorem \ref{thm 1.1}.

\hfill $\Box$ \\

\hspace{0.4cm}

\section{Limit of $\omega_{\epsilon }$-Hermitian-Einstein metrics }
\setcounter{equation}{0}

 Assume that  the reflexive Higgs sheaf $(\mathcal{E} , \phi )$ is $\omega $-stable. It is well known that the pulling back Higgs bundle $(E, \phi )$ is $\omega_{\epsilon }$-stable for sufficiently small $\epsilon$. By Simpson's result (\cite{Si}), there exists an $\omega_{\epsilon }$-Hermitian-Einstein metric $H_{\epsilon}$ for every  sufficiently small $\epsilon$.
 In this section,  we prove that, by choosing a subsequence and rescaling it, $H_{\epsilon}$ converges to an $\omega$-Hermitian-Einstein metric $H$ in local $C^{\infty}$-topology outside the exceptional divisor $\tilde{\Sigma }$.

 As above, let $\hat{H}$ be a fixed smooth Hermitian metric on the bundle $E$ over $\tilde{M}$. By taking a constant on $H_{\epsilon}$, we can suppose that
 \begin{equation}\label{det2}
  \int_{\tilde{M}}\tr \hat{S}_{\epsilon}\frac{\omega_{\epsilon}^{n}}{n!}=\int_{\tilde{M}}\log \det(\hat{h}_{\epsilon})\frac{\omega_{\epsilon}^{n}}{n!}=0.
 \end{equation}
 where $\exp (\hat{S}_{\epsilon})=\hat{h}_{\epsilon}=\hat{H}^{-1}H_{\epsilon }$.

 Let $H_{\epsilon } (t)$  be the long time solutions of  the heat flow (\ref{DDD1}) on the Higgs bundle $(E, \phi )$ with the fixed initial metric $\hat{H}$ and with respect to the K\"ahler metric $\omega_{\epsilon}$. We set:
   \begin{equation}\exp (\tilde{S}_{\epsilon}(t))=\tilde{h}_{\epsilon}(t)=H_{\epsilon}(t)^{-1}H_{\epsilon }.\end{equation}
 By (\ref{1}), (\ref{det2}) and noting that $\exp (\hat{S}_{\epsilon})=\exp (S_{\epsilon}(t))\exp (\tilde{S}_{\epsilon}(t))$, we have
 \begin{equation}\label{det3}
 \int_{\tilde{M}}\tr \tilde{S}_{\epsilon}(t)\frac{\omega_{\epsilon}^{n}}{n!}=\int_{\tilde{M}}\log \det(\tilde{h}_{\epsilon}(t))\frac{\omega_{\epsilon}^{n}}{n!}=0
 \end{equation}
  for all $ t\geq 0$. We first give a uniform $L^{1}$ estimate of $\hat{S}_{\epsilon }$.

\medskip

\begin{lemma}\label{lem 5.1}
There exists a constant $\hat{C}$ which is independent of $\epsilon$, such that
\begin{equation}\label{L102}
\|\hat{S}_{\epsilon }\|_{L^{1}(\tilde{M}, \omega_{\epsilon}, \hat{H})}:= \int_{\tilde{M}}|\hat{S}_{\epsilon }|_{\hat{H}}\frac{\omega_{\epsilon}^{n}}{n!}\leq \hat{C}
\end{equation}
for all $0<\epsilon \leq 1$.
\end{lemma}

\medskip

{\bf Proof. } We prove (\ref{L102}) by contradiction. If not, there exists a subsequence $\epsilon_{i}\rightarrow 0$ such that
\begin{equation}
\lim_{i\rightarrow \infty} \|\hat{S}_{\epsilon_{i} }\|_{L^{1}(\tilde{M}, \omega_{\epsilon_{i}}, \hat{H})}\rightarrow \infty .
\end{equation}
By (\ref{c02}), (\ref{C0a}) and (\ref{s01}), we also have
\begin{equation}
\lim_{i\rightarrow \infty} \|\tilde{S}_{\epsilon_{i} }(t)\|_{L^{1}(\tilde{M}, \omega_{\epsilon_{i}}, H_{\epsilon_{i}}(t))}\rightarrow \infty ,
\end{equation}
for all $t>0$.
 By (\ref{la02}), the uniform lower bound of Green functions $G_{\epsilon}$ (\ref{H005}) and the inequalities (\ref{c02}),
 we have
\begin{equation}\label{L103}
\|\tilde{S}_{\epsilon }(1)\|_{L^{\infty}(\tilde{M},  H_{\epsilon}(1))}\leq \grave{C}_{1}\|\tilde{S}_{\epsilon }(1)\|_{L^{1}(\tilde{M}, \omega_{\epsilon}, H_{\epsilon}(1))}+\grave{C}_{2},
\end{equation}
where $\grave{C}_{1}$ and $\grave{C}_{2}$ are uniform constants independent of $\epsilon$ and $t$. Using the inequality (\ref{s01}) again, we have
\begin{equation}\label{CM011}
\begin{split}
\|\tilde{S}_{\epsilon }(t )\|_{L^{\infty}(\tilde{M}, H_{\epsilon}(t))}\leq &\ r^{2}\grave{C}_{1}(\|\tilde{S}_{\epsilon }(t)\|_{L^{1}(\tilde{M}, \omega_{\epsilon}, H_{\epsilon}(t))}+\|S_{\epsilon }(t, 1)\|_{L^{1}(\tilde{M}, \omega_{\epsilon}, H_{\epsilon}(1))})\\ & +r\|S_{\epsilon }(t , 1)\|_{L^{\infty}(\tilde{M},  H_{\epsilon}(1))}+r\grave{C}_{2}\\
\end{split}
\end{equation}
 for all $ t >0$.

Set $\tilde{u}_{i}(t)=\|\tilde{S}_{\epsilon_{i} }(t)\|_{L^{1}(\tilde{M}, \omega_{\epsilon_{i}}, H_{\epsilon_{i}}(t))}^{-1}\tilde{S}_{\epsilon_{i} }(t)$, then $\|\tilde{u}_{i}(t)\|_{L^{1}(\tilde{M}, \omega_{\epsilon}, H_{\epsilon}(t))}=1$.
By (\ref{det3}) and (\ref{CM011}), we have $\int_{\tilde{M}}\tr u_{i}(t) \frac{\omega_{\epsilon}^{n}}{n!}=0$ and $\|\tilde{u}_{i}(t)\|_{L^{\infty}(\tilde{M}, H_{\epsilon_{i}}(t))}\leq C(t)$. Since $H_{\epsilon} (t) \rightarrow H(t)$ locally in $C^{\infty}$-topology and $\omega_{\epsilon }$ are locally uniform bounded outside $\tilde{\Sigma }$, by the Lemma 5.4 in \cite{Si}, we can show that, by choosing a subsequence which we also denote by $\tilde{u_i}(t)$, we have $\tilde{u}_{i}(t)\rightarrow \tilde{u}(t)$ weakly in $L_{1, loc}^{2}(\tilde{M}\setminus \tilde{\Sigma }, \omega, H(t) )$, where the limit $\tilde{u}(t)$ satisfies: $\|\tilde{u}(t)\|_{L^{1}(\tilde{M}\setminus \tilde{\Sigma }, \omega, H(t))}=1$, $\int_{\tilde{M}\setminus \tilde{\Sigma }}\tr(\tilde{u}(t))\frac{\omega^{n} }{n!}=0$. By (\ref{CM011}), we have \begin{equation}\label{tt01}
 \|\tilde{u}(t)\|_{L^{\infty}(\tilde{M}\setminus \tilde{\Sigma }, \omega, H(t))}\leq r^{2}\grave{C}_{1}.
 \end{equation}
Furthermore, if $\Upsilon : R\times R \rightarrow R$ is a positive smooth function such that $\Upsilon (\lambda_{1}, \lambda_{2})< (\lambda_{1}- \lambda_{2})^{-1}$ whenever $\lambda_{1}>\lambda_{2}$, then
 \begin{equation}\label{tt02}
 \begin{split}
 &\int_{\tilde{M}\setminus \tilde{\Sigma }}\tr (\tilde{u}(t)\sqrt{-1}\Lambda_{\omega }(F_{H(t), \phi })) + \langle\Upsilon  (\tilde{u}(t))(\overline{\partial }_{\phi }\tilde{u}(t)), \overline{\partial }_{\phi }\tilde{u}(t)\rangle_{H(t)}\frac{\omega^{n} }{n!}\\&\leq 0.\\
 \end{split}
 \end{equation}

Since $M\setminus \Sigma $ is biholomorphic to $\tilde{M}\setminus \tilde{\Sigma }$, and $\mathcal{E}$ is locally free on $M\setminus \Sigma $,  $\tilde{u}(t)$ can be seen as an $L_{1}^{2}$ section of $\End(\mathcal{E})$. By the same argument as that in section 4 (the proof of (\ref{semi02})), we can show that, by choosing a subsequence $t\rightarrow 0$, we have $\tilde{u}(t) \rightarrow \tilde{u}_{0}$ weakly in local $L_{1}^{2}$, where $\tilde{u}_{0}$ satisfies
\begin{equation}\int_{M}\tr(\tilde{u}_{0})\frac{\omega^{n} }{n!}=0, \quad  \|\tilde{u}_{0}\|_{L^{1}(M\setminus \Sigma , \omega , \hat{H})}=1,  \quad \|\tilde{u}(t)\|_{L^{\infty}(M\setminus \Sigma , \hat{H})}\leq r^{2}\grave{C}_{1}.\end{equation}
 \medskip
and
\begin{equation}\label{semi022}
 \int_{M\setminus \Sigma }\tr (\tilde{u}_{0}\sqrt{-1}\Lambda_{\omega }F_{\hat{H}, \phi }) + \langle\Upsilon(\tilde{u}_{0})(\overline{\partial }_{\phi}\tilde{u}_{0}), \overline{\partial }_{\phi}\tilde{u}_{0}\rangle_{\hat{H}}\frac{\omega^{n} }{n!}\leq 0.
 \end{equation}

 Now, by Simpson's trick (P888 in \cite{Si}), we can construct a saturated Higgs subsheaf $\mathcal{F}$ of $(\mathcal{E}, \phi )$ with $\mu_{\omega }(\mathcal{F})\geq \mu_{\omega }(\mathcal{E})$, which contradicts with the stability of  $(\mathcal{E}, \phi )$.

 \hfill $\Box$ \\

{\bf Proof of Theorem \ref{thm 1.2} } Since $\|\hat{S}_{\epsilon }\|_{L^{1}(\tilde{M},  \omega_{\epsilon}, \hat{M})}$ are uniformly bounded, by (\ref{c02}), (\ref{C0a}) and (\ref{s01}), there also exists a uniform constant $\grave{C}_{3}$ such that
\begin{equation}
\|\tilde{S}_{\epsilon }(1)\|_{L^{1}(\tilde{M}, \omega_{\epsilon}, H_{\epsilon}(1))}\leq \grave{C}_{3}.
\end{equation}
By (\ref{L103}), we have
\begin{equation}\label{L104}
\|\tilde{S}_{\epsilon }(1)\|_{L^{\infty}(\tilde{M},  H_{\epsilon}(1))}\leq \grave{C}_{1}\grave{C}_{3}+\grave{C}_{2}
\end{equation}
 for all $0< \epsilon \leq 1$. By the local estimate (\ref{C01}) in Lemma \ref{lem 2.3}, we see that there exists a constant $\tilde{C}_{0}(\delta^{-1})$ independent of $\epsilon$ such that
 \begin{equation}
|\hat{S}_{\epsilon }|_{\hat{H}}(x)\leq \tilde{C}_{0}(\delta^{-1})
\end{equation}
for all $x\in \tilde{M}\setminus B_{\omega_1}(\delta)$ and all $0<\epsilon \leq 1$. Since $H_{\epsilon}$ satisfies the $\omega_{\epsilon}$-Hermitian-Einstein equation (\ref{HE}), by the same argument as that in Lemmas \ref{lem 2.4} and \ref{lem 2.5} in section 2, we have  uniform higher-order estimates for $h_{\epsilon }$, i.e. there exist constants $\tilde{C}_{k}(\delta^{-1})$ independent of $\epsilon$, such that
\begin{equation}
\|\hat{h}_{\epsilon}\|_{C^{k+1,\alpha}, \tilde{M}\setminus B_{\omega_1}(2\delta)}\leq \tilde{C}_{k+1}(\delta^{-1})
\end{equation}
for all $k\geq 0$ and all $0<\epsilon \leq 1$. So by choosing a subsequence, we have $H_{\epsilon} $
converges to a Hermitian metric $H $ on $M\setminus \Sigma $ in locally $C^{\infty}$-topology, and $H$ satisfies the Hermitian-Einstein equation, i.e.
\begin{equation}
\sqrt{-1}\Lambda_{\omega }(F_{H}+[\phi , \phi ^{\ast H}])=\lambda \Id_{\mathcal{E}}.
\end{equation}

By (\ref{L104}), we see that the metrics $H(1)$ and $H$  are  mutually bounded each other on $\mathcal{E}|_{M\setminus \Sigma}$. On the other hand,  we have shown that $|\phi|_{H(1), \omega}\in L^{\infty}(M)$ in section 3, then $|\phi |_{H, \omega}$  also belongs to $L^{\infty}(M)$. This implies that $|\Lambda_{\omega }(F_{H})|_{H}$ is uniform bounded on $M\setminus \Sigma $. By (\ref{CW22}), it is easy to see that $|F_{H}|_{H, \omega}$ is square integrable. So we know that the metric $H$ is an admissible Hermitian-Einstein metric on the Higgs sheaf $(\mathcal{E}, \phi )$. This completes the proof of Theorem \ref{thm 1.2}.

\hfill $\Box$ \\

\hspace{0.3cm}

\end{document}